\PassOptionsToPackage{unicode}{hyperref}
\PassOptionsToPackage{hyphens}{url}
\PassOptionsToPackage{dvipsnames,svgnames,x11names}{xcolor}
\documentclass[
]{article}
\usepackage{xcolor}
\usepackage{amsmath,amssymb}
\usepackage{iftex}
\ifPDFTeX
  \usepackage[T1]{fontenc}
  \usepackage[utf8]{inputenc}
  \usepackage{textcomp} % provide euro and other symbols
\else % if luatex or xetex
  \usepackage{unicode-math} % this also loads fontspec
  \defaultfontfeatures{Scale=MatchLowercase}
  \defaultfontfeatures[\rmfamily]{Ligatures=TeX,Scale=1}
\fi
 \usepackage{lmodern}
 \ifPDFTeX\else
 \fi
\IfFileExists{upquote.sty}{\usepackage{upquote}}{}
\IfFileExists{microtype.sty}{% use microtype if available
  \usepackage[]{microtype}
  \UseMicrotypeSet[protrusion]{basicmath} % disable protrusion for tt fonts
}{}
\makeatletter
\@ifundefined{KOMAClassName}{% if non-KOMA class
  \IfFileExists{parskip.sty}{%
    \usepackage{parskip}
  }{% else
    \setlength{\parindent}{0pt}
    \setlength{\parskip}{6pt plus 2pt minus 1pt}}
}{% if KOMA class
  \KOMAoptions{parskip=half}}
\makeatother
\makeatletter
\ifx\paragraph\undefined\else
  \let\oldparagraph\paragraph
  \renewcommand{\paragraph}{
    \@ifstar
      \xxxParagraphStar
      \xxxParagraphNoStar
  }
  \newcommand{\xxxParagraphStar}[1]{\oldparagraph*{#1}\mbox{}}
  \newcommand{\xxxParagraphNoStar}[1]{\oldparagraph{#1}\mbox{}}
\fi
\ifx\subparagraph\undefined\else
  \let\oldsubparagraph\subparagraph
  \renewcommand{\subparagraph}{
    \@ifstar
      \xxxSubParagraphStar
      \xxxSubParagraphNoStar
  }
  \newcommand{\xxxSubParagraphStar}[1]{\oldsubparagraph*{#1}\mbox{}}
  \newcommand{\xxxSubParagraphNoStar}[1]{\oldsubparagraph{#1}\mbox{}}
\fi
\makeatother

\usepackage{color}
\usepackage{fancyvrb}

\DefineVerbatimEnvironment{Highlighting}{Verbatim}{commandchars=\\\{\}}
\usepackage{framed}
\definecolor{shadecolor}{RGB}{241,243,245}
\newenvironment{Shaded}{\begin{snugshade}}{\end{snugshade}}

\newcommand{\BuiltInTok}[1]{\textcolor[rgb]{0.00,0.23,0.31}{#1}}

\newcommand{\DecValTok}[1]{\textcolor[rgb]{0.68,0.00,0.00}{#1}}

\newcommand{\FloatTok}[1]{\textcolor[rgb]{0.68,0.00,0.00}{#1}}

\newcommand{\ImportTok}[1]{\textcolor[rgb]{0.00,0.46,0.62}{#1}}

\newcommand{\NormalTok}[1]{\textcolor[rgb]{0.00,0.23,0.31}{#1}}
\newcommand{\OperatorTok}[1]{\textcolor[rgb]{0.37,0.37,0.37}{#1}}

\newcommand{\SpecialCharTok}[1]{\textcolor[rgb]{0.37,0.37,0.37}{#1}}
\newcommand{\SpecialStringTok}[1]{\textcolor[rgb]{0.13,0.47,0.30}{#1}}

\usepackage{longtable,booktabs,array}
\usepackage{calc} % for calculating minipage widths
\usepackage{etoolbox}
\makeatletter
\patchcmd\longtable{\par}{\if@noskipsec\mbox{}\fi\par}{}{}
\makeatother
\IfFileExists{footnotehyper.sty}{\usepackage{footnotehyper}}{\usepackage{footnote}}
\makesavenoteenv{longtable}
\usepackage{graphicx}
\makeatletter
\newsavebox\pandoc@box
\newcommand*\pandocbounded[1]{% scales image to fit in text height/width
  \sbox\pandoc@box{#1}%
  \Gscale@div\@tempa{\textheight}{\dimexpr\ht\pandoc@box+\dp\pandoc@box\relax}%
  \Gscale@div\@tempb{\linewidth}{\wd\pandoc@box}%
  \ifdim\@tempb\p@<\@tempa\p@\let\@tempa\@tempb\fi% select the smaller of both
  \ifdim\@tempa\p@<\p@\scalebox{\@tempa}{\usebox\pandoc@box}%
  \else\usebox{\pandoc@box}%
  \fi%
}
\def\fps@figure{htbp}
\makeatother

\NewDocumentCommand\citeproctext{}{}

\makeatletter
 \let\@cite@ofmt\@firstofone
 \def\@biblabel#1{}
 \def\@cite#1#2{{#1\if@tempswa , #2\fi}}
\makeatother
\newlength{\cslhangindent}
\newlength{\csllabelwidth}
\newenvironment{CSLReferences}[2] % #1 hanging-indent, #2 entry-spacing
 {\begin{list}{}{%
  \setlength{\itemindent}{0pt}
  \setlength{\leftmargin}{0pt}
  \setlength{\parsep}{0pt}
  \ifodd #1
   \setlength{\leftmargin}{\cslhangindent}
   \setlength{\itemindent}{-1\cslhangindent}
  \fi
  \setlength{\itemsep}{#2\baselineskip}}}
 {\end{list}}
\usepackage{calc}

\providecommand{\tightlist}{%
  \setlength{\itemsep}{0pt}\setlength{\parskip}{0pt}}

\usepackage{arxiv}
\usepackage{orcidlink}
\usepackage{amsmath}
\usepackage[T1]{fontenc}
\makeatletter
\@ifpackageloaded{caption}{}{\usepackage{caption}}
\AtBeginDocument{%
\ifdefined\contentsname
  \renewcommand*\contentsname{Table of contents}
\else
  \newcommand\contentsname{Table of contents}
\fi
\ifdefined\listfigurename
  \renewcommand*\listfigurename{List of Figures}
\else
  \newcommand\listfigurename{List of Figures}
\fi
\ifdefined\listtablename
  \renewcommand*\listtablename{List of Tables}
\else
  \newcommand\listtablename{List of Tables}
\fi
\ifdefined\figurename
  \renewcommand*\figurename{Figure}
\else
  \newcommand\figurename{Figure}
\fi
\ifdefined\tablename
  \renewcommand*\tablename{Table}
\else
  \newcommand\tablename{Table}
\fi
}
\@ifpackageloaded{float}{}{\usepackage{float}}
\floatstyle{ruled}
\@ifundefined{c@chapter}{\newfloat{codelisting}{h}{lop}}{\newfloat{codelisting}{h}{lop}[chapter]}
\floatname{codelisting}{Listing}

\usepackage{amsthm}
\theoremstyle{definition}
\newtheorem{definition}{Definition}[section]
\theoremstyle{plain}
\newtheorem{lemma}{Lemma}[section]
\theoremstyle{plain}
\newtheorem{theorem}{Theorem}[section]
\theoremstyle{definition}
\newtheorem{example}{Example}[section]
\theoremstyle{remark}
\AtBeginDocument{}

\newtheorem{refremark}{Remark}[section]

\makeatother
\makeatletter
\@ifpackageloaded{caption}{}{\usepackage{caption}}
\@ifpackageloaded{subcaption}{}{\usepackage{subcaption}}
\makeatother
\usepackage{bookmark}
\IfFileExists{xurl.sty}{\usepackage{xurl}}{} % add URL line breaks if available
\hypersetup{
  pdftitle={Multi-Objective Optimization with Desirability and Morris-Mitchell Criterion},
  pdfkeywords={Heuristic algorithms, Multi-objective continuous
optimization, Software for MCDM, Applications -
Engineering, Desirability function, Design of experiments, Space-filling
design, Multi-objective optimization, Surrogate
modeling, Morris-Mitchell criterion, Maximin criterion, Minimax
criterion, Sequential parameter optimization, Infill-point plots},
  colorlinks=true,
  linkcolor={blue},
  filecolor={Maroon},
  citecolor={Blue},
  urlcolor={Blue},
  pdfcreator={LaTeX via pandoc}}

\newcommand{\runninghead}{A Preprint (Version 2)}
\title{Multi-Objective Optimization with Desirability and
Morris-Mitchell Criterion}
\usepackage{etoolbox}
\makeatletter
\providecommand{\subtitle}[1]{% add subtitle to \maketitle
  \apptocmd{\@title}{\par {\large #1 \par}}{}{}
}
\makeatother
\subtitle{An Experimental Analysis}
\def\asep{\\\\\\ } % default: all authors on same column
\def\asep{\And }
\author{\textbf{Thomas
Bartz-Beielstein}~\orcidlink{0000-0002-5938-5158}\\\\Bartz \& Bartz
GmbH, 51643 Gummersbach,
Germany\\\\\href{mailto:bartzbeielstein@gmail.com}{bartzbeielstein@gmail.com}\asep\textbf{Eva
Bartz}\\\\Bartz \& Bartz GmbH, 51643 Gummersbach,
Germany\\\\\href{mailto:eva.bartz@bartzundbartz.de}{eva.bartz@bartzundbartz.de}\asep\textbf{Alexander
Hinterleitner}\\\\Bartz \& Bartz GmbH, 51643 Gummersbach,
Germany\\\\\href{mailto:ahinterleitner@bartzundbartz.de}{ahinterleitner@bartzundbartz.de}\asep\textbf{Christoph
Leitenmeier}\\\\Everllence SE, Engineering Turbocharger, 86153 Augsburg,
Germany\\\\\href{mailto:christoph.leitenmeier@everllence.com}{christoph.leitenmeier@everllence.com}\asep\textbf{Ihab
Abd El Hussein}\\\\Everllence SE, Engineering Turbocharger, 86153
Augsburg,
Germany\\\\\href{mailto:ihab.abd-el-hussein@everllence.com}{ihab.abd-el-hussein@everllence.com}}
\date{}
\begin{document}
\maketitle
\begin{abstract}
Industrial experimental designs frequently lack optimal space-filling
properties, rendering them unrepresentative. This study presents a
comprehensive methodology to refine existing designs by enhancing
coverage quality while optimizing experimental outcomes. We discuss and
analyse variants of the Morris-Mitchell criterion to quantify and
improve spatial distributions. Based on potential theory, we analyze
monotonicity properties and limitations of the Morris-Mitchell criteria.
Practically, we implement a multi-objective optimization framework
utilizing the Python packages \texttt{spotdesirability} and
\texttt{spotoptim}. This framework uses desirability functions to
combine surrogate-model predictions with space-filling enhancements into
a unified score. Demonstrated through data from a compressor development
case study, this approach optimizes performance objectives alongside
design coverage. To facilitate implementation, we introduce novel
infill-point diagnostics that visually guide the sequential placement of
design points. This integrated methodology successfully bridges spatial
theory with engineering application, balancing the crucial exploration
and exploitation trade-off.
\end{abstract}
{\bfseries \emph Keywords}
\def\sep{\textbullet\ }
Heuristic algorithms \sep Multi-objective continuous
optimization \sep Software for MCDM \sep Applications -
Engineering \sep Desirability function \sep Design of
experiments \sep Space-filling design \sep Multi-objective
optimization \sep Surrogate modeling \sep Morris-Mitchell
criterion \sep Maximin criterion \sep Minimax criterion \sep Sequential
parameter optimization \sep 
Infill-point plots

\section{Introduction}\label{introduction}

For planned experiments, the $n$ design points, forming an
experimental design $X_p$, can be chosen to optimally cover the input
space\footnote{ In the following, $X_p$ denotes the planned design,
  and $X_u$ denotes the unplanned (non-optimal) design, which is
  already available.}. Depending on the model used, there are various
approaches for this. Classical approaches are the optimality criteria
used in Design of Experiments (DoE), which were developed for linear
models (A-optimality, D-optimality, E-optimality,
T-optimality)\footnote{ A-optimality minimizes the average variance of
  the parameter estimates (trace of the inverse information matrix).
  D-optimality maximizes the determinant of the information matrix,
  minimizing the volume of the confidence ellipsoid for the parameters.
  E-optimality maximizes the minimum eigenvalue of the information
  matrix, minimizing the worst-case parameter estimation variance.
  T-optimality maximizes the trace of the information matrix, focusing
  on maximizing the total information about all parameters.}. In the
field of Design and Analysis of Computer Experiments (DACE), where
models such as Kriging are used, space-filling designs have become
established (Santner et al. 2003).

However, these optimality criteria are not directly applicable if
already existing, non-optimal designs should be used and extended. In
practice, results from unplanned, arbitrarily chosen designs, say
$X_u$, that do not satisfy any optimality criteria are often already
available. These designs are biased, i.e., they are primarily based on
the experience of the experimenters and usually cover only a small part
of the input space and therefore do not meet the statistical optimality
criteria. As a consequence, the results of the experiments are not
representative. The question therefore arises as to how these designs
can be improved to increase the quality of the coverage of the input
space while simultaneously improving the results of the experiments.
This leads to a multiobjective optimization problem.

We present an approach for extending existing, non-optimal designs,
which is based on the Morris-Mitchell criterion $\Phi$ (Morris and
Mitchell 1995). This criterion measures the distances between the design
points in the input space to quantify the quality of the coverage of the
input space. The criterion can be used to distribute the design points
uniformly in the input space. However, $\Phi$ grows with the number of
design points $n$, so that it is not directly applicable to existing
designs. To address this issue, we use the intensified version of the
criterion, denoted by $\Phi^{\ast}$. It takes into account the number
of design points (Bartz-Beielstein 2025a). We demonstrate how
$\Phi^{\ast}$ can be used to improve the quality of the coverage of
the input space. A theoretical analysis based on potential theory shows
that $\Phi^{\ast}$ is a more appropriate criterion for existing
designs, but still has some limitations. Therefore, we also consider as
a second option a corrected version of the Morris-Mitchell criterion,
denoted by $\hat{\Phi}$.

To demonstrate the applicability of our approach, we use a case study
based on a problem from compressor development. In designs from
industry, datasets (i.e., data rows = number of investigated designs)
are always limited, yet statements should be made about all possible
designs within the investigated input parameter space. The particular
challenge is that individual input/output parameters often exhibit
pronounced local minimum/maximum effects, so it is usually unclear
whether this local design space has been discretized finely enough.

We define a multi-objective desirability function that combines the
predictions from multiple surrogate models and $\Phi^{\ast}$ into a
single desirability score. The problem used in the case study is an
already multi-objective optimization problem with nine different
objectives. In practice, modeling using desirability functions is a
common method to optimize multiple objectives simultaneously (Derringer
and Suich 1980; National Institute of Standards and Technology 2021). It
was therefore natural to combine the Morris-Mitchell criterion
$\Phi^{\ast}$ with desirability functions to quantify the quality of
the coverage of the input space and to distribute the design points
uniformly in the input space as suggested by Bartz-Beielstein (2025a).
Thus, a total of ten objectives are available (nine multiobjective ones
from the application plus one for the space-filling property) that can
be used in the case study. We limit ourselves here to the two objectives
that have the greatest relevance in practice and add $\Phi^{\ast}$ as
an additional objective. This article demonstrates how to perform this
multi-objective optimization task using \texttt{spotoptim}. To visualize
the impact of a new design point on existing design, we will use
infill-point plots (Bartz-Beielstein 2025a). Infill-point plots
visualize the location of new design points in the input space
relatively to existing design points.

This article is structured as follows: Section~\ref{sec-dataset}
describes the dataset used in the case study. The surrogate models are
introduced in Section~\ref{sec-surrogate-models}. The Morris-Mitchell
criterion and its variants are explained in
Section~\ref{sec-MM-criterion}. Desirability functions are described in
Section~\ref{sec-desirability}. Section~\ref{sec-optimization} presents
the multi-objective optimization approach. Section~\ref{sec-opt-wo-mm}
shows the optimization without the Morris-Mitchell criterion, whereas
Section~\ref{sec-mm-criterion-multi} describes the optimization with the
Morris-Mitchell criterion. This article concludes with presentation of
the results and a discussion in Section~\ref{sec-results-discussion}.

\section{The Dataset from Compressor Development}\label{sec-dataset}

\subsection{Data Loading and
Preparation}\label{data-loading-and-preparation}

The project is based on dimensionless characteristic numbers from
various numerically executed compressor designs, which were collected in
an industrial context.

The independent (feature) variables are stored in the pandas DataFrame
\texttt{df\_x\_normalized}. The dependent (target) variables are stored
in the pandas DataFrame \texttt{df\_z\_normalized}. These data are
freely available and can be loaded as follows:

\protect\phantomsection\label{load-data-spotdesirability}
\begin{Shaded}
\begin{Highlighting}[]
\ImportTok{import}\NormalTok{ spotdesirability.data\_utils }\ImportTok{as}\NormalTok{ du}
\NormalTok{df\_x\_normalized, df\_z\_normalized }\OperatorTok{=}\NormalTok{ du.load\_compressor\_data()}
\end{Highlighting}
\end{Shaded}

For the nine target variables, the desirabilities for the optimal values
shown in Table~\ref{tbl-desirability} are used. More than two hundred
designs were available, exactly $n=213$ as tabulated data. Here, each
row corresponds to a unique design, and each of the $k=27$ columns to
an individual variable. We consider these 27 independent input variables
or features (equivalent to geometric parameters such as dimensions and
angles) and nine target variables ($p=9$) or outputs (mainly resulting
performance parameters). All data, i.e., input and target variables, are
normalized to a range of $[0,1]$.

Two target variables are selected for the multi-objective optimization
problem discussed in this article: the 8th and the first column (index 7
and 0, i.e., column \texttt{z8} and \texttt{z1}) of the target data
frame \texttt{y}. Therefore, we will use \texttt{x1}, \texttt{x2},
\ldots, \texttt{x27} as feature names and \texttt{z8}, and \texttt{z1}
as target names. These desirabilities are derived from practical
application and reflect the requirements for the target variables. Since
the data are provided in the accompanied dataset, the reader can perform
similar case studies with different target variables and desirabilities.

\begin{longtable}[]{@{}llrc@{}}
\caption{The nine target variables and their
desirabilities}\label{tbl-desirability}\tabularnewline
\toprule\noalign{}
Index & Label & Desirability & Used \\
\midrule\noalign{}
\endfirsthead
\toprule\noalign{}
Index & Label & Desirability & Used \\
\midrule\noalign{}
\endhead
\bottomrule\noalign{}
\endlastfoot
0 & z1 & max & yes \\
1 & z2 & max & \\
2 & z3 & max & \\
3 & z4 & max & \\
4 & z5 & max & \\
5 & z6 & 0.5 & \\
6 & z7 & min & \\
7 & z8 & max & yes \\
8 & z9 & min & \\
\end{longtable}

\begin{figure}

\centering{

\pandocbounded{\includegraphics[keepaspectratio]{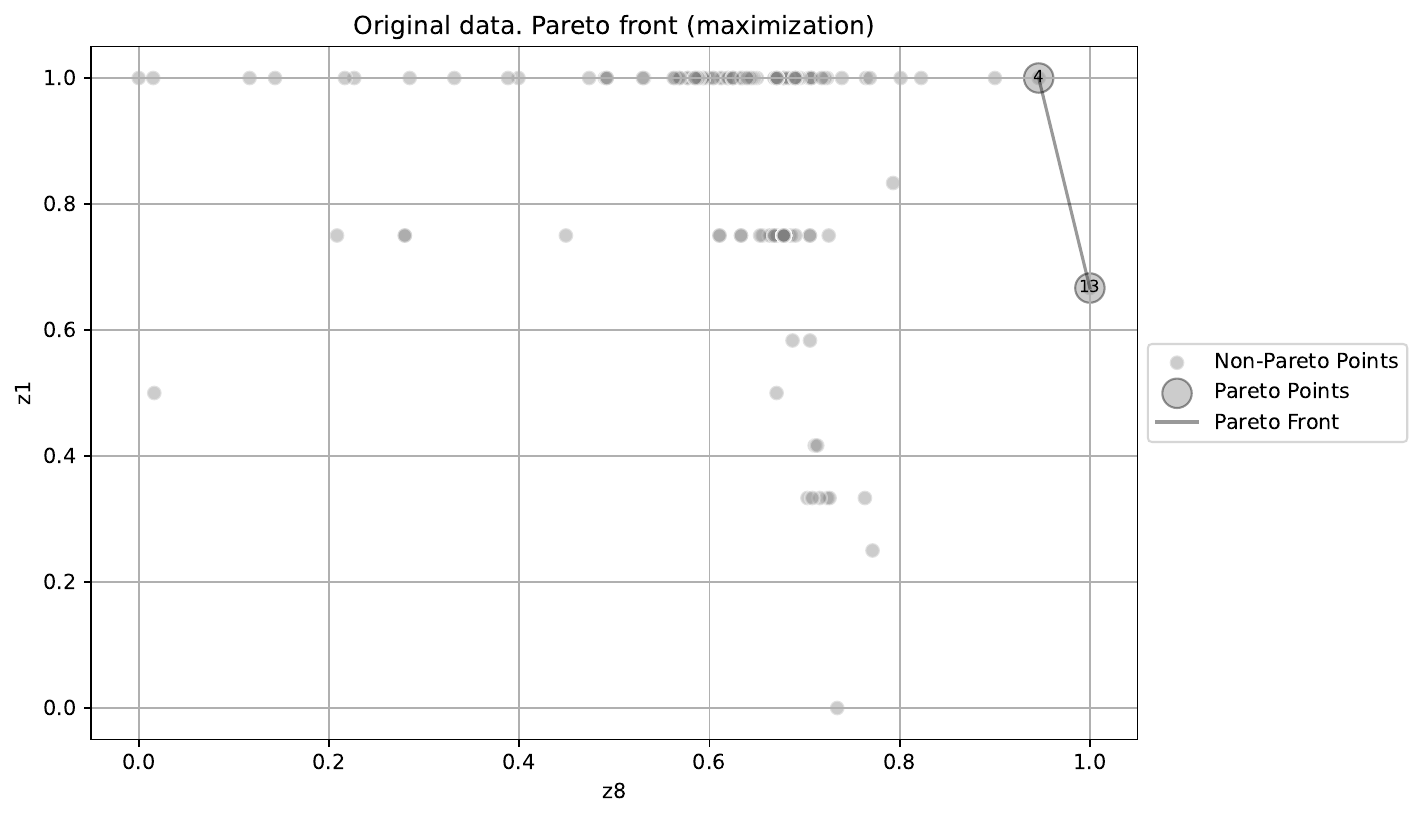}}

}

\caption{\label{fig-paretofront-orig-ap3-x}Pareto front of the original
data. Both functions should be maximized}

\end{figure}

The visualization of the Pareto front of the original data is performed
using the function \texttt{plot\_mo} from the software package
\texttt{spotoptim}. This function creates a 2D visualization of the
Pareto fronts for the target variable pairs. When viewing the Pareto
front (Figure~\ref{fig-paretofront-orig-ap3-x}), points \texttt{4} and
\texttt{13} stand out due to high \texttt{z8} values and could be
classified as outliers. Since these are real-world data, stemming from
existing designs, it is not clear whether these points are truly
outliers or not. Thus, we did not remove these points. But further
investigation of this data may be useful.

\section{Surrogate Models}\label{sec-surrogate-models}

To accomplish the goal of this study, we will use surrogate models to
predict the target variables. Based on these predictions, we will
perform multi-objective optimization to determine the next infill point.
This point will be proposed to the experimenter, who will evaluate it
and provide the corresponding target values. Before performing these
costly experiments, the information from the surrogate models will be
used to estimate its quality. New in our approach is the use of
desirability functions to estimate the quality of the next infill point
in combination with the intensified Morris-Mitchell criterion. We have
selected two surrogate models for this study: a Random Forest Regressor
and a Gaussian Process Regressor.

\subsection{Random Forest Regressor}\label{random-forest-regressor}

We use a Random Forest Regressor from \texttt{scikit-learn} with default
hyper-parameters. Important hyper-parameters are the number of trees in
the forest (\texttt{n\_estimators}), the maximum depth of the trees
(\texttt{max\_depth}), and the random state for reproducibility
(\texttt{random\_state}). Random Forest Regressor models from
\texttt{sckit-learn} can natively handle multiple-target problems, see
the discussion in
\url{https://scikit-learn.org/stable/auto_examples/ensemble/plot_random_forest_regression_multioutput.html}.
So, we can alteratively train one single model for the two target
variables \texttt{z8} and \texttt{z1}, or one multi-objective model for
both variables. Since results are identical, we have chosen the latter
approach to guarantee compatibility with other algorithms which do not
allow multi-target computations.

\subsection{The Evaluation Function for Multiple
Models}\label{the-evaluation-function-for-multiple-models}

The function \texttt{mo\_eval\_models()} from the \texttt{spotoptim}
package is used to evaluate multiple models for multi-output regression.
The data set is randomly split into training and test sets, where
\texttt{X\_train} and \texttt{y\_train} are used for training the
models, and \texttt{X\_test} and \texttt{y\_test} are used for
evaluation. The \texttt{test\_size} parameter defines the proportion of
the data used for testing. It is set to \texttt{0.3} here.

\begin{figure}
\centering{
\pandocbounded{\includegraphics[width=0.8\linewidth,keepaspectratio]{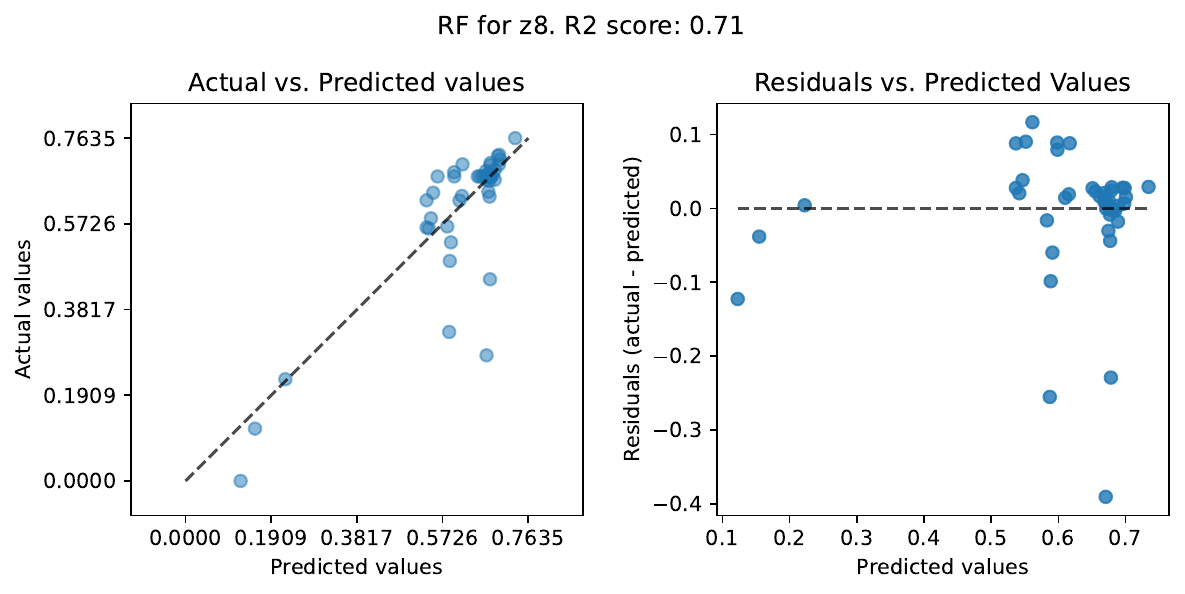}}
}
\caption{\label{fig-plot_predictions-z8-multi-rf}Predictions for z8 of
the multi-objective RF model compared to the original data}
\end{figure}

\begin{figure}
\centering{
\pandocbounded{\includegraphics[width=0.8\linewidth,keepaspectratio]{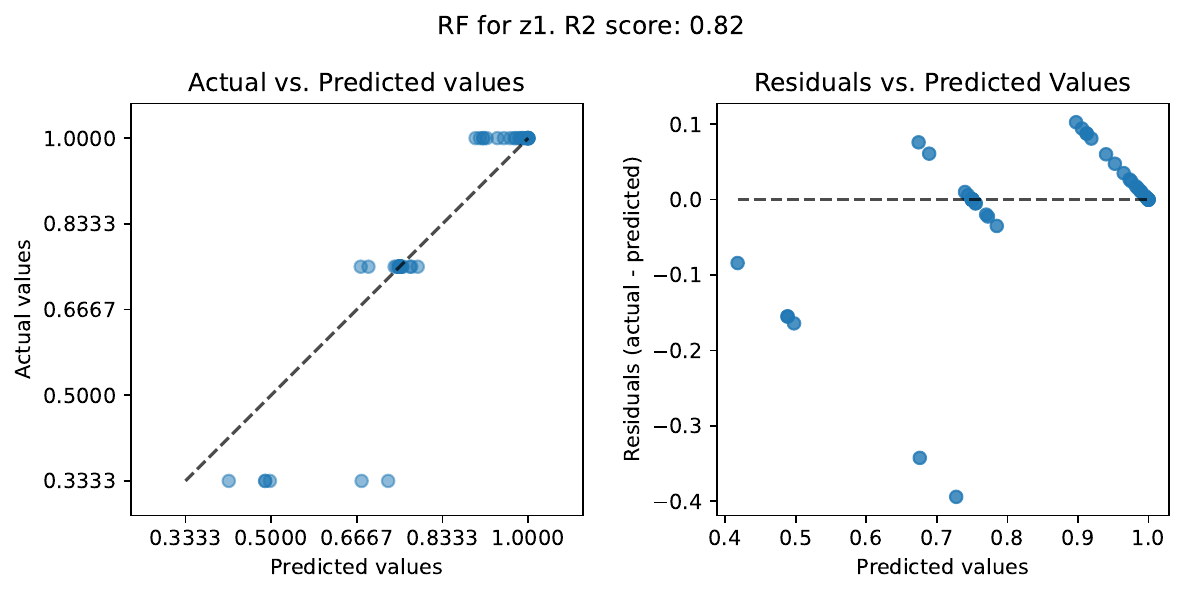}}
}
\caption{\label{fig-plot_predictions-z1-multi-rf}Predictions for z1 of
the multi-objective RF model compared to the original data}

\end{figure}

\begin{figure}
\centering{
\pandocbounded{\includegraphics[width=\linewidth,keepaspectratio]{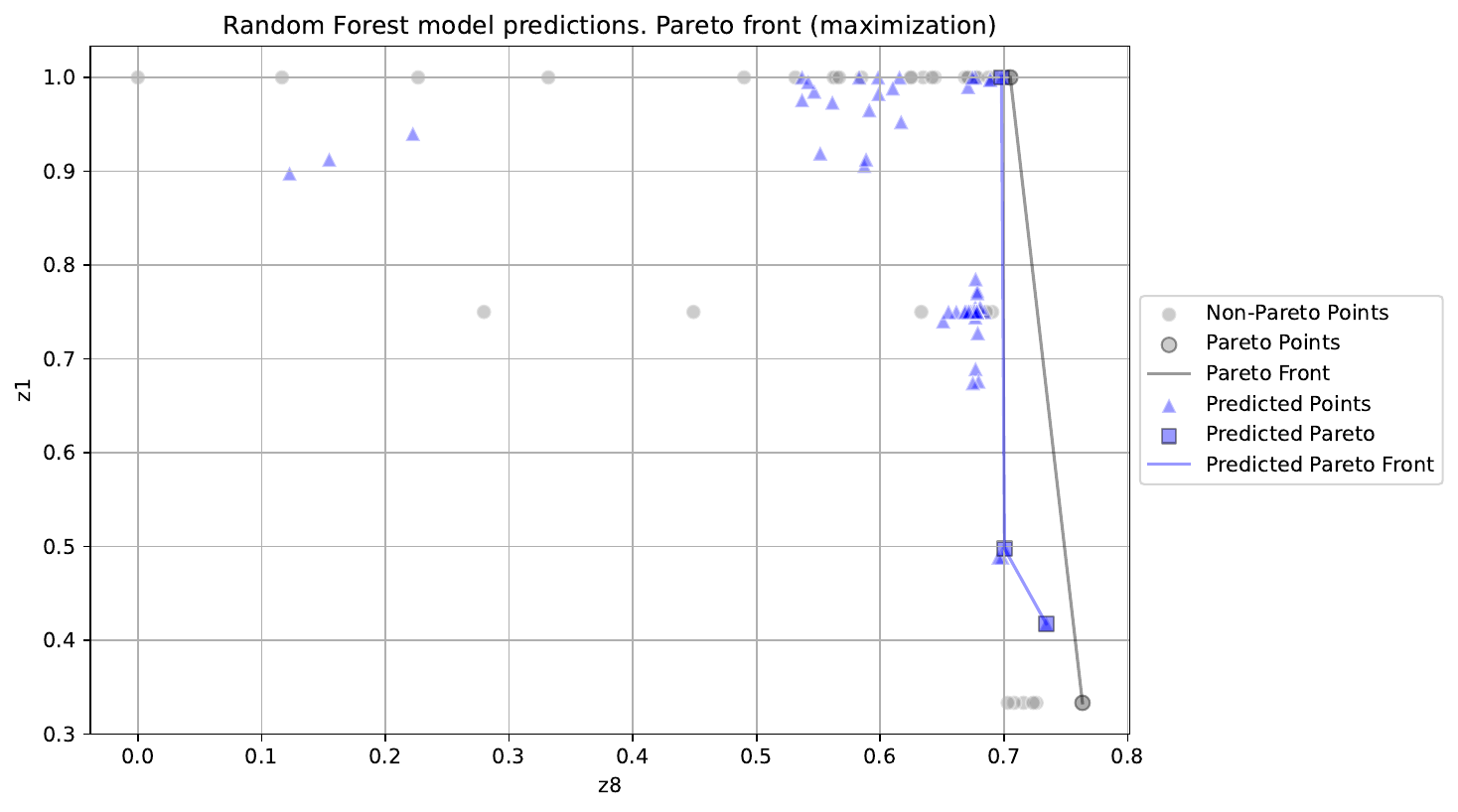}}
}
\caption{\label{fig-paretofront-rf-ap3-x-plot_mo}Pareto front of the
Random Forest surrogate model predictions and the original Pareto front
(with possible outliers)}

\end{figure}

Figure~\ref{fig-plot_predictions-z8-multi-rf} and
Figure~\ref{fig-plot_predictions-z1-multi-rf} show the predictions for
the target variables \texttt{z8} and \texttt{z1} of the random forest
surrogate model compared to the original data. The Pareto fronts of the
predictions and the original data are shown in
Figure~\ref{fig-paretofront-rf-ap3-x-plot_mo}. From
Figure~\ref{fig-paretofront-gp-ap3-x-plot_mo} we can conclude that the
surrogate-model based approach is able to ``rediscover'' the original
Pareto front.

\subsection{Gaussian Process
Regressor}\label{gaussian-process-regressor}

\begin{figure}

\centering{

\pandocbounded{\includegraphics[keepaspectratio]{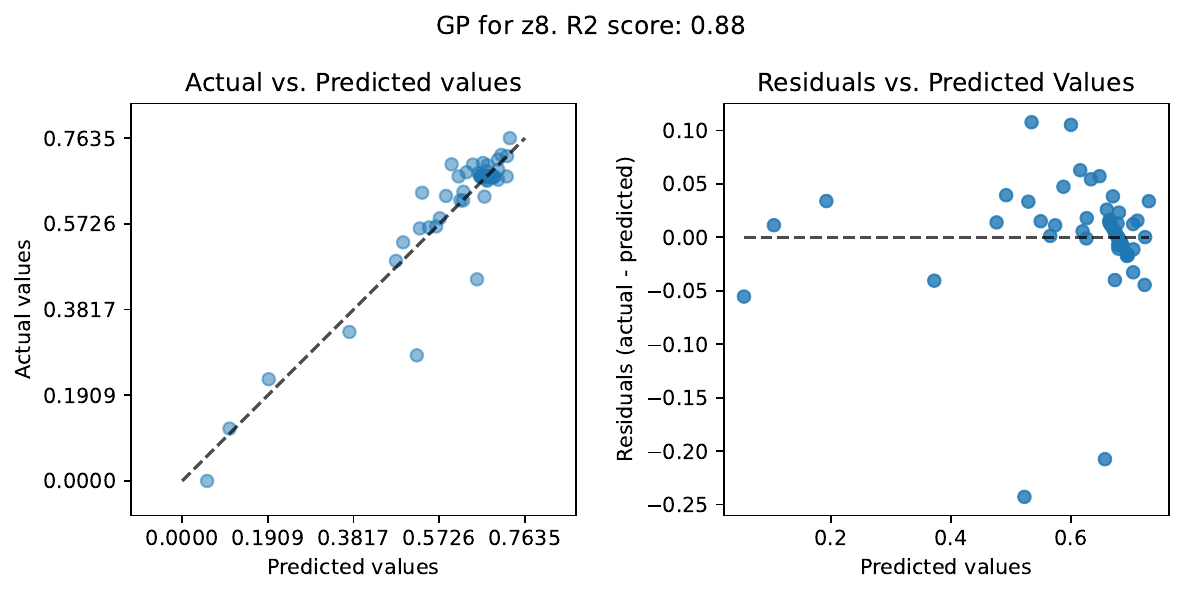}}

}

\caption{\label{fig-plot_predictions-z8-multi-gp}Predictions for z8 of
the multi-objective GP model compared to the original data}

\end{figure}

\begin{figure}
\centering{
\pandocbounded{\includegraphics[width=0.8\linewidth,keepaspectratio]{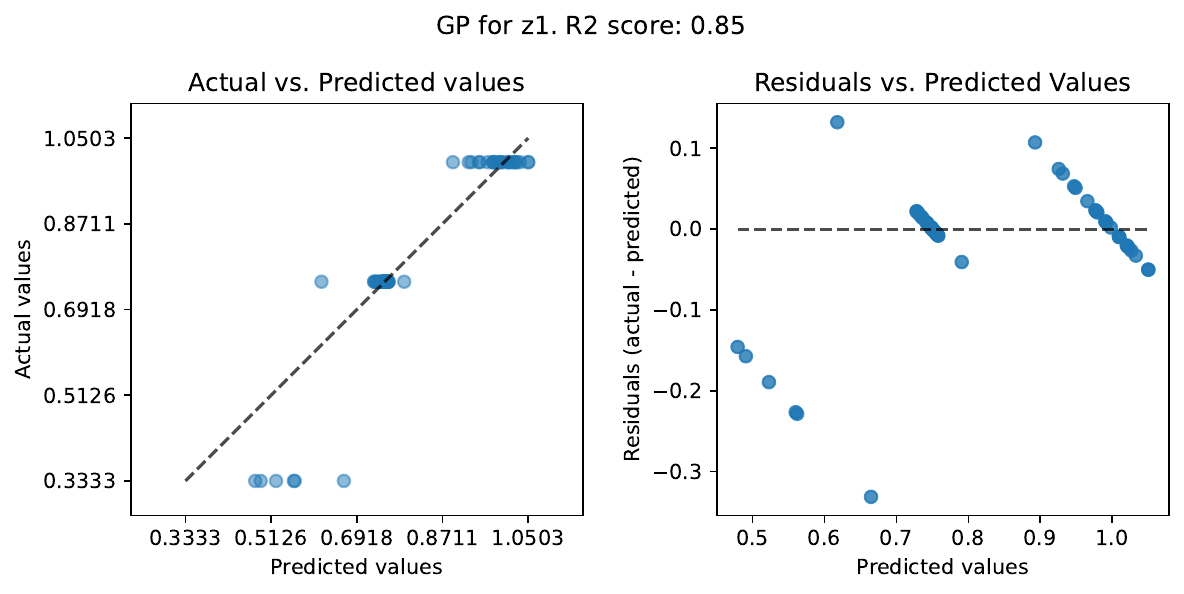}}
}
\caption{\label{fig-plot_predictions-z1-multi-gp}Predictions for z1 of
the multi-objective GP model compared to the original data}

\end{figure}

To enable multi-output regression for Gaussian Process (GP) models,
similar to Random Forest, we fit a separate GP (optionally with Nystroem
kernel approximation\footnote{The Nystroem method works by randomly
  selecting \texttt{n\_components} training samples as basis points,
  computing the kernel matrix on that subset, and using it to produce an
  approximate feature map with \texttt{n\_components} output dimensions.})
for each target variable and stack the predictions. Similar to the
visual analysis of the Random Forest model, we can plot the predictions
for both target variables \texttt{z8} and \texttt{z1}, which are shown
in Figure~\ref{fig-plot_predictions-z8-multi-gp} and
Figure~\ref{fig-plot_predictions-z1-multi-gp}, respectively.

\begin{figure}

\centering{

\pandocbounded{\includegraphics[width=\linewidth,keepaspectratio]{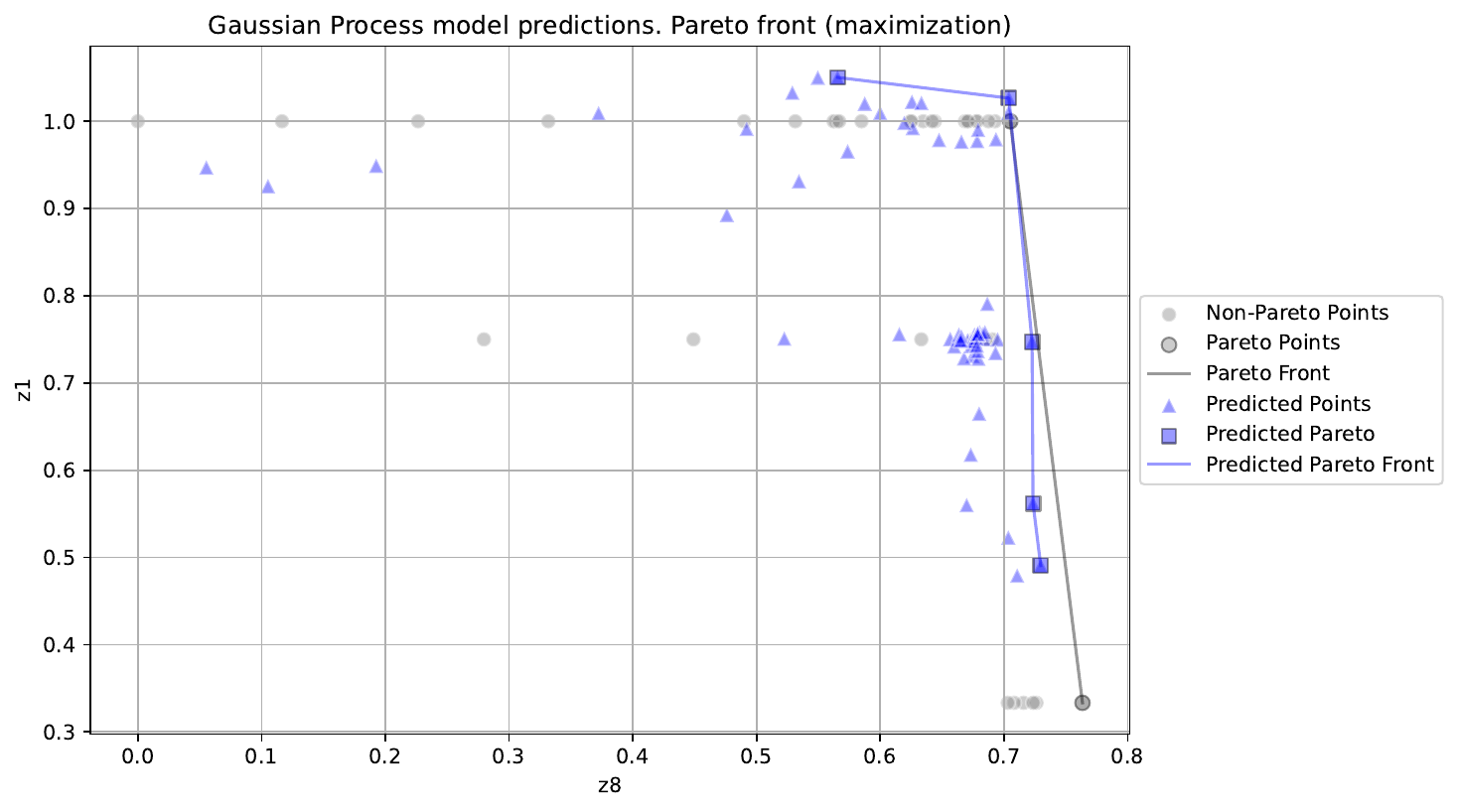}}

}

\caption{\label{fig-paretofront-gp-ap3-x-plot_mo}Pareto front of the
Gaussian Process model predictions}

\end{figure}

The Pareto fronts of the predictions of the Gaussian Process model and
the original data are displayed in
Figure~\ref{fig-paretofront-gp-ap3-x-plot_mo}.

\subsection{Cross-validation Comparison of the
Models}\label{cross-validation-comparison-of-the-models}

Because the evaluation on the train-test-splitted data is not a good
indicator for the generalization performance of the models, a
cross-validated comparison of the Random Forest and Gaussian Process
models is shown in the appendix (Section~\ref{sec-cv-model-comparison}).
The cross-validated comparison shows that the Random Forest model
outperforms the Gaussian Process model and generates more robust
solutions. In addition, the Random Forest model is also much faster to
train. Therefore, we decided to use the Random Forest model for the
further analysis.

\section{An Intensified Morris-Mitchell
Criterion}\label{sec-MM-criterion}

The Morris-Mitchell criterion describes the quality of an experimental
design based on the distances between sample points. It aims to maximize
the minimum distance between any two points in the design (Morris and
Mitchell 1995). This ensures that no two points are too close to each
other, which promotes better coverage. We consider the default
Morris-Mitchell criterion $\Phi$ and first, before we introduce the
intensified version $\Phi^{\ast}$ and the concentrated version
$\hat{\Phi}$, which take into account the number of design points.

\subsection{The Default Morris-Mitchell
Criterion}\label{the-default-morris-mitchell-criterion}

The Morris-Mitchell criterion is a widely used criterion for evaluating
the space-filling properties of a sampling plan (also known as a design
of experiments).

\begin{definition}[Morris-Mitchell
Criterion]\protect\hypertarget{def-mmphi}{}\label{def-mmphi}

The Morris-Mitchell criterion, $\Phi$, is defined as:

\begin{equation}\protect\phantomsection\label{eq-mmphi}{
\Phi = \Phi_q(n) = \left( \sum_{i=1}^{m} J_i d_i^{-q} \right)^{1/q}
}\end{equation}

where:

\begin{itemize}
\tightlist
\item
  $n$ is the number of points in the sampling plan.
\item
  $d_1 < d_2 < \dots < d_m$ are the distinct distances between all
  pairs of points in the sampling plan.
\item
  $J_i$ is the number of pairs of points separated by the distance
  $d_i$.
\item
  $q$ is a large positive integer (e.g., $q=2, 5, 10, \dots$).
\end{itemize}

\end{definition}

As $q \to \infty$, minimizing $\Phi_q$ is equivalent to maximizing
the minimum distance $d_1$, and then minimizing the number of pairs
$J_1$ at that distance, and so on.

The Morris-Mitchell criterion is a maximin criterion.

\begin{definition}[Maximin and Minimax
Criteria]\protect\hypertarget{def-maximin-minimax}{}\label{def-maximin-minimax}

Maximin is the principle of designing experiments to maximize the
minimum pairwise distance between design points, i.e., 
$\max_{X} \min_{i \neq j} \, d(x_i, x_j)$.
As $q \to \infty$, minimizing $\Phi_q$ is equivalent to this --- it
pushes all points apart so that the smallest inter-point distance is as
large as possible. This prevents clustering. Minimax (also called
coverage criterion) is the dual concept: minimize the maximum distance
from any point in the design space to the nearest sample point, i.e.,
$\min_{X} \max_{x \in \mathcal{X}} \min_{i} \, d(x, x_i)$.
\end{definition}

Minimax criteria ensure no region of the input space is left uncovered.
So, while maximin looks inward at distances \emph{between} design
points; minimax looks outward at distances from the design to the rest
of the space. Both promote uniformity, but they are not equivalent: a
good maximin design is not necessarily a good minimax design and vice
versa.

The Morris-Mitchell criterion refines the pure maximin criterion by also
penalizing designs that have many point pairs at the minimum distance
(via the weights $J_1$), making it a richer characterization of
space-filling quality. The Morris-Mitchell criterion is widely used in
various fields, including global sensitivity analysis, uncertainty
quantification, and experimental design. Small values of the
Morris-Mitchell criterion indicate a better design. We use the
implementation of the Morris-Mitchell criterion from the
\texttt{spotoptim} package, which provides the function
\texttt{mmphi()}, which is based on the code provided by Forrester et
al. (2008).

\begin{figure}
\centering{
\pandocbounded{\includegraphics[keepaspectratio]{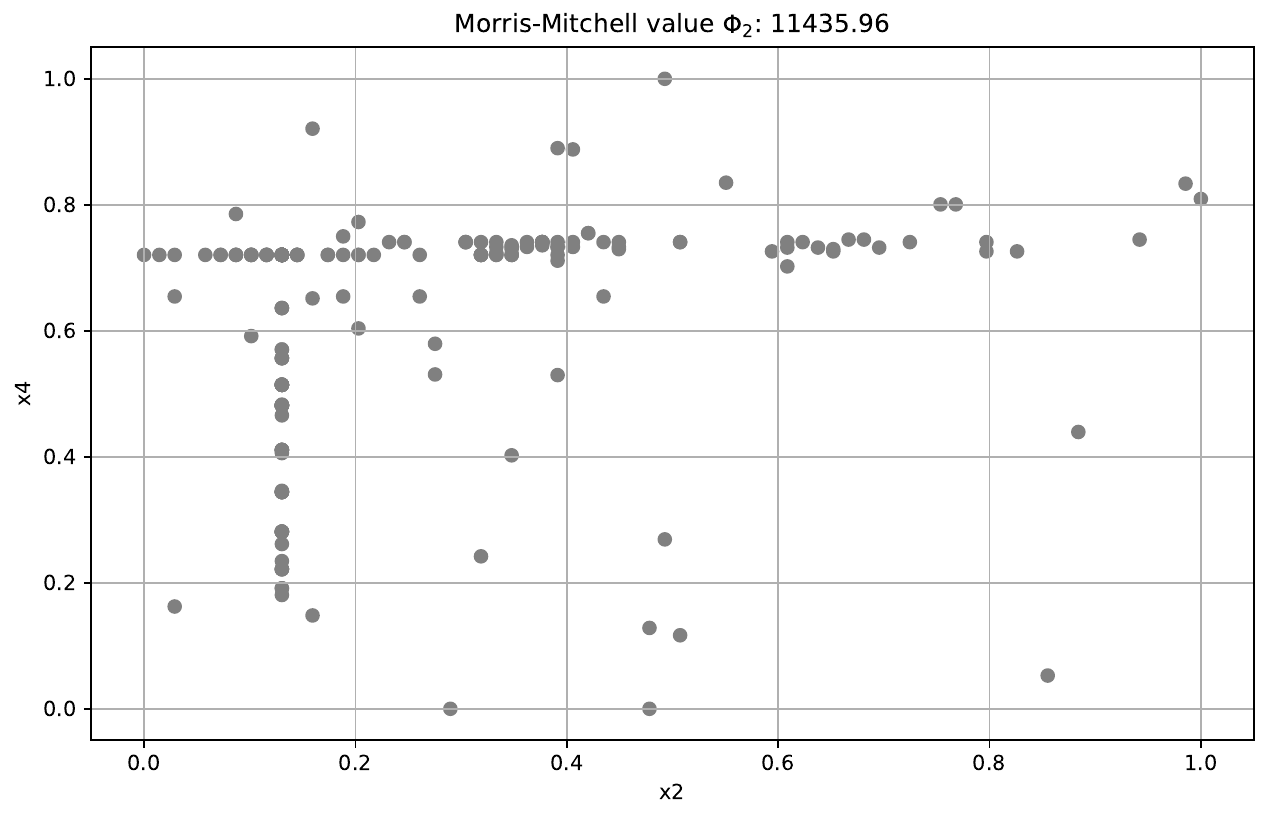}}
}
\caption{\label{fig-man-design}Design of Industrial Data. Typical
clusters and ``lanes'' are visible. These data have a relatively high
Morris-Mitchell value.}
\end{figure}

\begin{figure}
\centering{
\pandocbounded{\includegraphics[keepaspectratio]{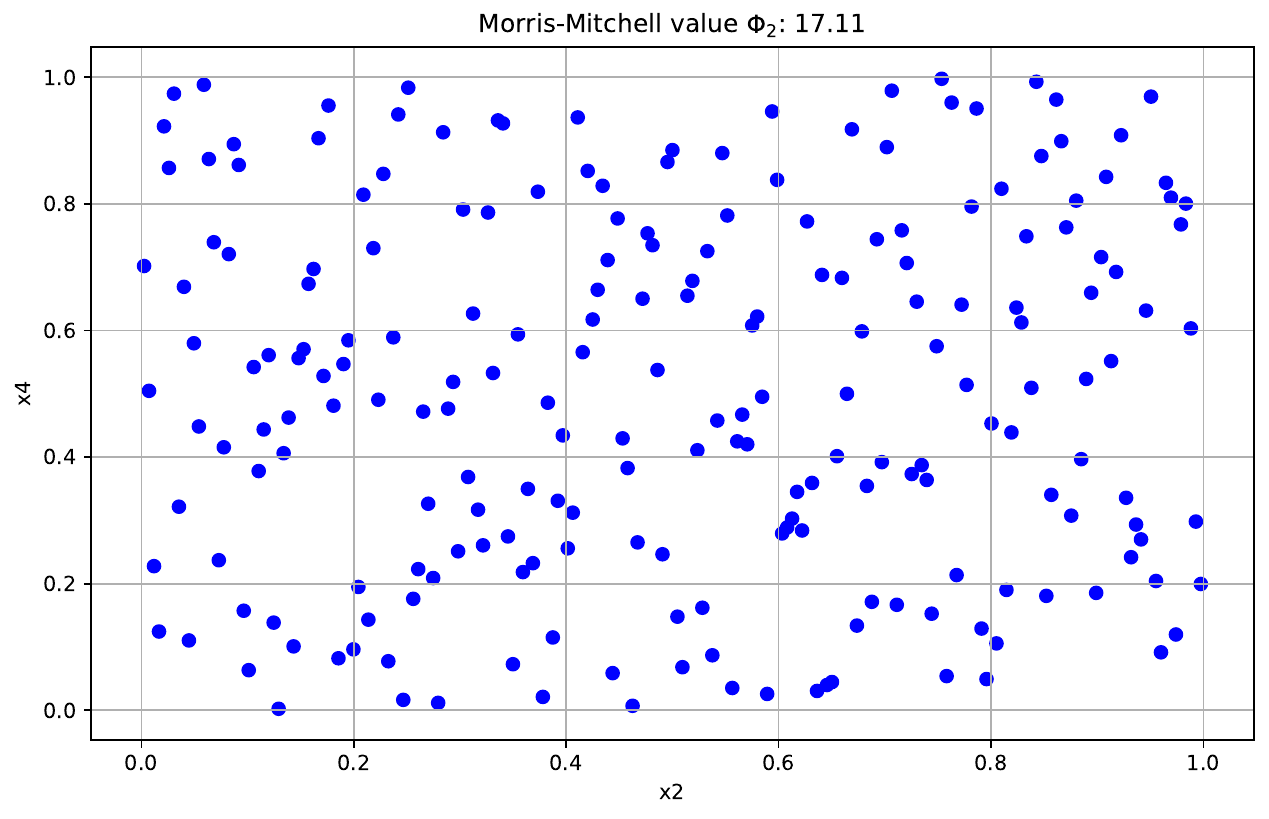}}
}
\caption{\label{fig-mm-design}Optimized design. Same design size as in
the industrial data.}
\end{figure}

Figure~\ref{fig-man-design} visualizes a design of industrial data with
n=213 sample points. It shows a scatter plot of the normalized input
variables \texttt{x2} and \texttt{x4}. This visualization was presented
by the industrial partner to illustrate the problem of non-space
fillingness. Experimenters proceed with great caution to avoid machine
failure or unfeasible settings that result in unpredictable system
behavior. The design does not exhibit good space-filling properties, as
there are clusters of points and large empty regions. For comparison,
Figure~\ref{fig-mm-design} shows an optimized space-filling design with
the same number of points (n=213). This design was computed using the
Morris-Mitchell criterion, which is used as an objective function in an
optimization algorithm. The \texttt{spotoptim}'s \texttt{bestlh()}
function searches for a sampling design that maximizes the criterion,
resulting in a well-distributed set of sample points.

\subsection{The Intensified Morris-Mitchell
Criterion}\label{the-intensified-morris-mitchell-criterion}

One limitation of the standard $\Phi_q$ is that its value depends on
$n$, the number of points in the sampling plan. This makes it
difficult to compare the quality of sampling plans with different sample
sizes. To address this, we can use an intensified version of the
criterion.

\begin{definition}[Intensified Morris-Mitchell
Criterion]\protect\hypertarget{def-mmphi-intensive}{}\label{def-mmphi-intensive}

The intensified criterion is defined as:

\begin{equation}\protect\phantomsection\label{eq-mmphi-intensive}{
\Phi^{\ast} = \Phi_q^{\ast}(n) = \left( \frac{1}{M} \sum_{i=1}^{m} J_i d_i^{-q} \right)^{1/q},
}\end{equation}

with $M = \binom{n}{2} = \frac{n(n-1)}{2}$, the total number of pairs
of points in the design.

\end{definition}

This intensification (or normalization) aims at allowing a fairer
comparison between designs of different sizes. A lower value should
indicate a better space-filling property. But this is not always the
case, as will be shown in the following.

\subsection{\texorpdfstring{Using \texttt{mmphi\_intensive} in
\texttt{spotoptim}}{Using mmphi\_intensive in spotoptim}}\label{using-mmphi_intensive-in-spotoptim}

The \texttt{spotoptim} library provides the \texttt{mmphi\_intensive}
function to calculate this criterion.

\begin{example}[Usage of
\texttt{mmphi\_intensive}]\protect\hypertarget{exm-mmphi-intensive}{}\label{exm-mmphi-intensive}

We create a simple three-point sampling plan in 2D (n=3) and calculate
the intensified Morris-Mitchell criterion with q=2, using Euclidean
distances (p=2).

\end{example}

\protect\phantomsection\label{mmphi-intensive-example}
\begin{Shaded}
\begin{Highlighting}[]
\NormalTok{X3 }\OperatorTok{=}\NormalTok{ np.array([}
\NormalTok{    [}\FloatTok{0.0}\NormalTok{, }\FloatTok{0.0}\NormalTok{],}
\NormalTok{    [}\FloatTok{0.5}\NormalTok{, }\FloatTok{0.5}\NormalTok{],}
\NormalTok{    [}\FloatTok{1.0}\NormalTok{, }\FloatTok{1.0}\NormalTok{]}
\NormalTok{])}
\NormalTok{quality, J, d }\OperatorTok{=}\NormalTok{ mmphi\_intensive(X3, q}\OperatorTok{=}\DecValTok{2}\NormalTok{, p}\OperatorTok{=}\DecValTok{2}\NormalTok{)}
\BuiltInTok{print}\NormalTok{(}\SpecialStringTok{f"Quality (Phi\_q\_intensive): }\SpecialCharTok{\{}\NormalTok{quality}\SpecialCharTok{\}}\SpecialStringTok{"}\NormalTok{)}
\BuiltInTok{print}\NormalTok{(}\SpecialStringTok{f"Multiplicities (J): }\SpecialCharTok{\{}\NormalTok{J}\SpecialCharTok{\}}\SpecialStringTok{"}\NormalTok{)}
\BuiltInTok{print}\NormalTok{(}\SpecialStringTok{f"Distinct Distances (d): }\SpecialCharTok{\{}\NormalTok{d}\SpecialCharTok{\}}\SpecialStringTok{"}\NormalTok{)}
\end{Highlighting}
\end{Shaded}

\begin{verbatim}
Quality (Phi_q_intensive): 1.224744871391589
Multiplicities (J): [2 1]
Distinct Distances (d): [0.70710678 1.41421356]
\end{verbatim}

The returned \texttt{quality} is the $\Phi_q^{\ast}$ value. The array
\texttt{J} contains the multiplicities for each distinct distance, and
\texttt{d} contains the distinct distances found in the design. There
are three points, resulting in three distances, $d_{12}$, $d_{13}$,
and $d_{23}$. The smallest distance is $d_{12} = d_{23}$, and it
occurs twice, so $J_1 = 2$. The largest distance is $d_{13}$, and it
occurs once, so $J_2 = 1$.

\subsubsection{\texorpdfstring{Efficient Updates with
\texttt{mmphi\_intensive\_update}}{Efficient Updates with mmphi\_intensive\_update}}\label{efficient-updates-with-mmphi_intensive_update}

When constructing a design sequentially (e.g., adding one point at a
time), recalculating the full distance matrix and criterion from scratch
can be inefficient. The \texttt{mmphi\_intensive\_update} function from
the \texttt{spotoptim} package allows updating $\Phi_q^{\ast}$
efficiently by only computing distances between the new point and the
existing points.

\begin{example}[Usage of
\texttt{mmphi\_intensive\_update}]\protect\hypertarget{exm-mmphi-intensive-update}{}\label{exm-mmphi-intensive-update}

Assume we have the state from the previous example: X, quality, J, d,
and we want to add a new point.

\end{example}

\protect\phantomsection\label{mmphi-intensive-update-example-onepoint}
\begin{Shaded}
\begin{Highlighting}[]
\NormalTok{new\_point }\OperatorTok{=}\NormalTok{ np.array([}\FloatTok{0.1}\NormalTok{, }\FloatTok{0.1}\NormalTok{])}
\NormalTok{new\_quality, new\_J, new\_d }\OperatorTok{=}\NormalTok{ mmphi\_intensive\_update(X3, new\_point, J, d, q}\OperatorTok{=}\DecValTok{2}\NormalTok{, p}\OperatorTok{=}\DecValTok{2}\NormalTok{)}
\BuiltInTok{print}\NormalTok{(}\SpecialStringTok{f"Updated Quality: }\SpecialCharTok{\{}\NormalTok{new\_quality}\SpecialCharTok{\}}\SpecialStringTok{"}\NormalTok{)}
\end{Highlighting}
\end{Shaded}

\begin{verbatim}
Updated Quality: 3.115613474919968
\end{verbatim}

The function \texttt{mmphi\_intensive\_update} returns the updated
criterion, multiplicities, and distinct distances, which can be used for
the next update step.

\begin{figure}

\centering{

\pandocbounded{\includegraphics[keepaspectratio]{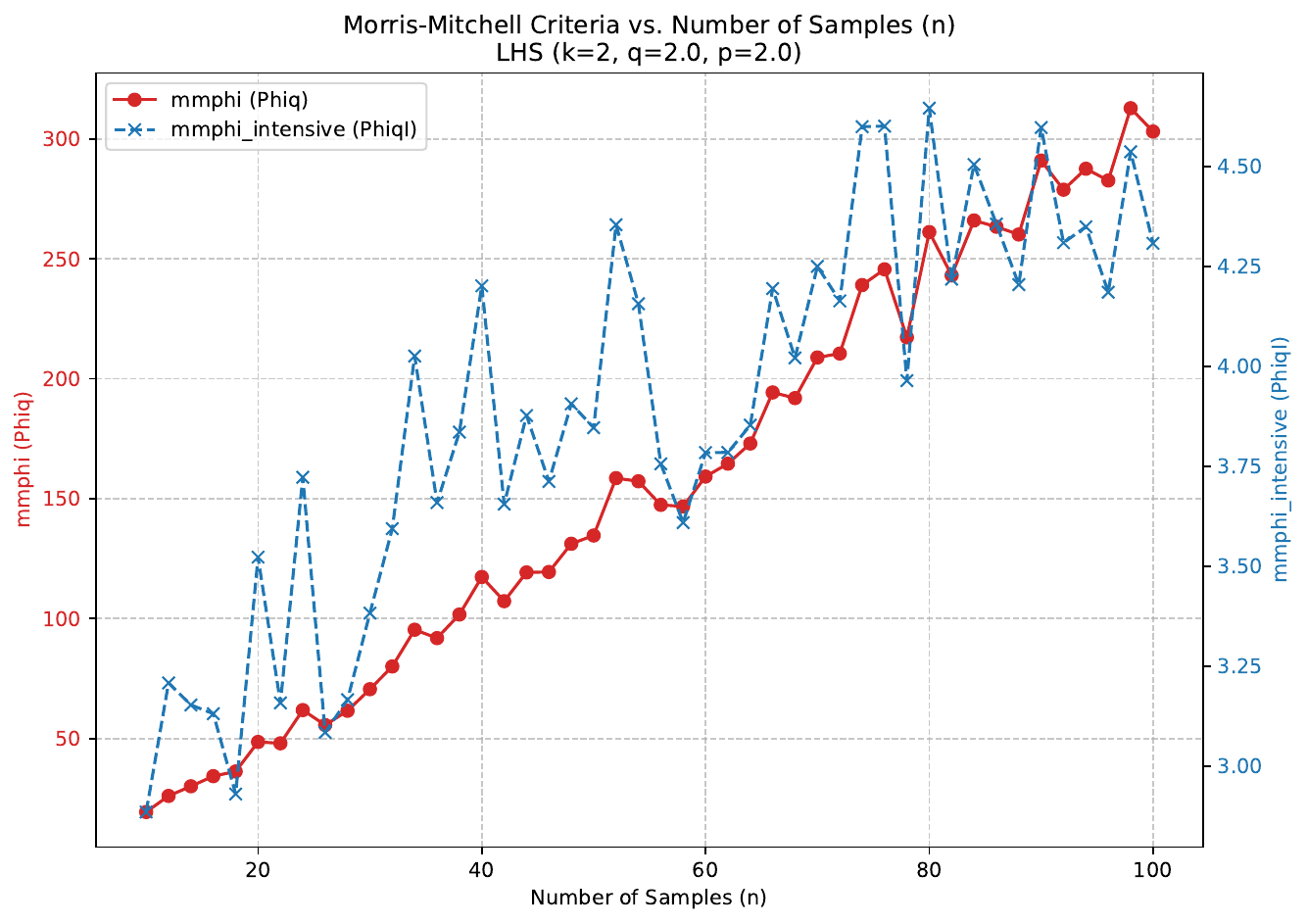}}

}

\caption{\label{fig-mm-fig-mm-vs-n-lhs}Two-dimensional LHD designs.
Morris-Mitchell criterion and extended criterion vs.~number of samples
(n). Note the different scales of the y-axes.}

\end{figure}

Figure~\ref{fig-mm-fig-mm-vs-n-lhs} shows the results of the
Morris-Mitchell criterion and its intensified version versus the number
of samples ($n$) for Latin Hypercube Sampling (LHS) designs in a
2-dimensional space. It was generated with the function
\texttt{plot\_mmphi\_vs\_n\_lhs} from the \texttt{spotoptim} package.
Figure~\ref{fig-mm-fig-mm-vs-n-lhs} shows that $\Phi_q^{\ast}$ is less
sensitive to the number of samples. The former raises from 50 to 300,
while the extended criterion raises only from 2.5 to 4.5.

\subsubsection{Sensitivity of the Morris-Mitchell Criterion with Respect
to Added Points}\label{sec-mm-sensitivity}

\begin{figure}
\centering{
\pandocbounded{\includegraphics[width=0.8\linewidth,keepaspectratio]{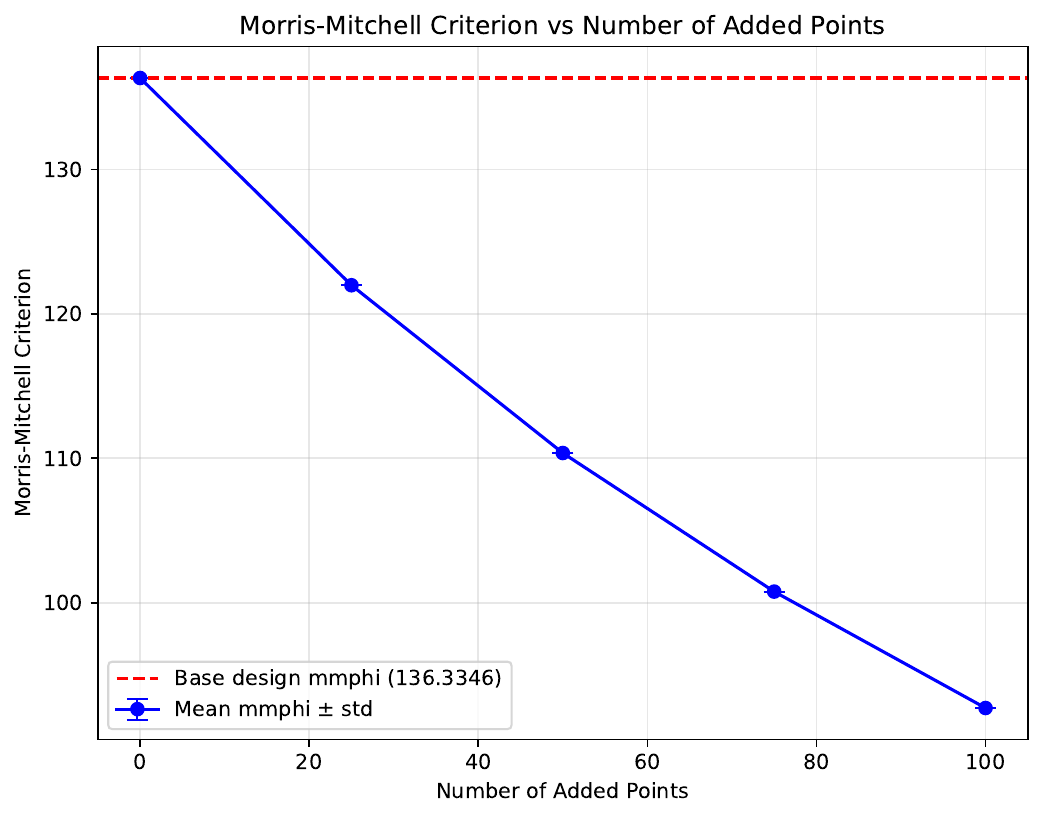}}
}
\caption{\label{fig-mmphi-intens-plot-x-call}Compressor data set. The
intensified Morris-Mitchell criterion $\Phi_q^\ast$ vs.~number of
added points. Randomized points are added to the entire dataset, which
results in a decrease of the criterion.}

\end{figure}

First, the quality of the current design $X_u$, which contains the
entire compressor dataset, is calculated using the intensified
Morris-Mitchell criterion $\Phi^{\ast}$. As shown in
Figure~\ref{fig-mmphi-intens-plot-x-call}, which shows $\Phi^{\ast}$
versus the number of added points,the initial design has a value of
$\Phi^{\ast} \approx 136$, which is used as a reference for the
optimization. Randomized points are added to the entire dataset. By
adding ten random points, the Morris-Mitchell criterion decreases from
$\Phi^{\ast} \approx 136$ to $\Phi^{\ast} \approx 130$. Adding one
hundred random points, the criterion further decreases to values below
$\Phi^{\ast} \approx 100$.

\begin{figure}

\centering{

\pandocbounded{\includegraphics[width=0.8\linewidth,keepaspectratio]{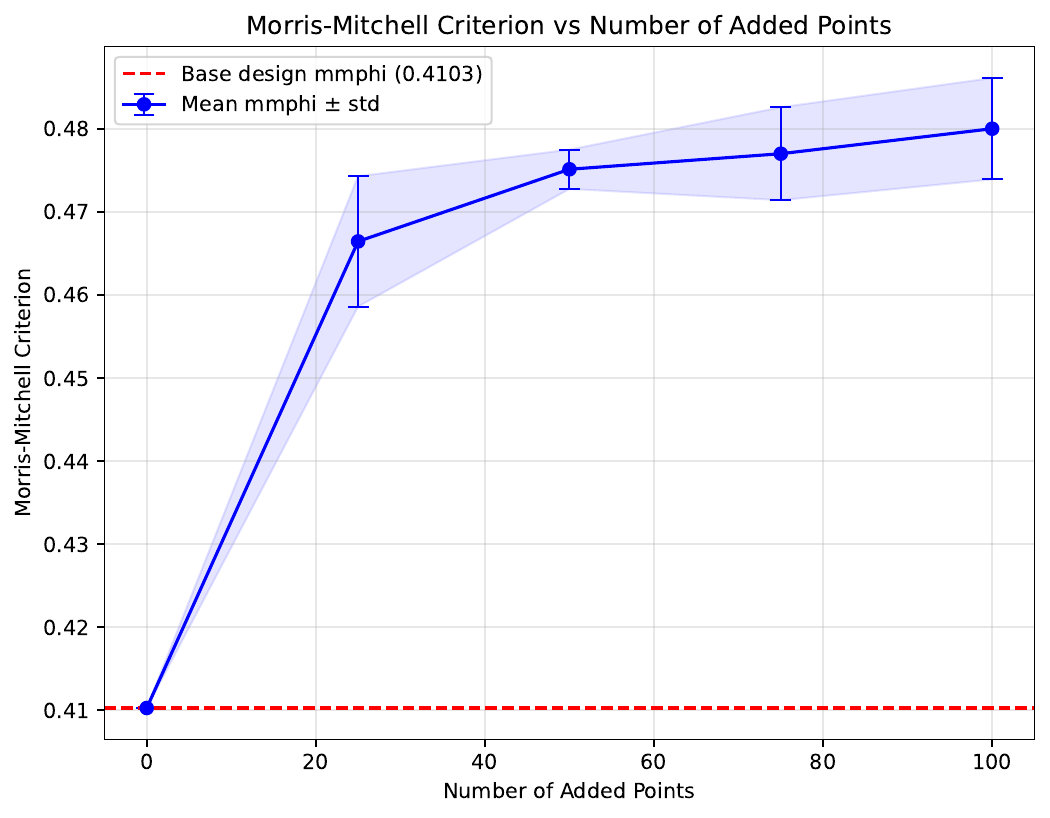}}

}

\caption{\label{fig-mmphi-intens-plot-x-search-scaled1}Optimal design.
The intensified Morris-Mitchell criterion vs.~number of added points.
Randomized points are added to an LHD, which results in an increase of
the criterion.}

\end{figure}

In contrast to the previous example, $\Phi^{\ast}$ is calculated for a
planned and optimized experimental design, i.e., random points are added
to an existing LHD $X_p$, which is nearly optimal.
Figure~\ref{fig-mmphi-intens-plot-x-search-scaled1} shows the
$\Phi^{\ast}$ values after adding random points to the existing design
$X_p$. Here, $\Phi^{\ast}$ increases from $\approx 0.41$ to
$\approx 0.43$ after adding ten random points.

These experiments motivated the usage of the Morris-Mitchell criterion
to improve the quality of the coverage of the input space, which can be
formulated as Remark~\ref{rem-mmphi-sensitivity}.

\begin{refremark}[Hypothesis on the Sensitivity of the Morris-Mitchell
Criterion to Added Points]
\leavevmode

\begin{itemize}
\tightlist
\item
  If a ``good point'' (w.r.t. the design criterion) is added to a ``bad
  design'', the Morris-Mitchell value can be reduced, so that the
  difference between the $\Phi^{\ast}$ value of the current design and
  the $\Phi^{\ast}$ value of the design with the added point is
  positive and should be maximized. This was illustrated in
  Figure~\ref{fig-mmphi-intens-plot-x-call}.
\item
  If a ``bad point'' (w.r.t. the design criterion) is added to a ``good
  design'', the $\Phi^{\ast}$ value can slightly increase, as
  illustrated in Figure~\ref{fig-mmphi-intens-plot-x-search-scaled1}.
\end{itemize}

\label{rem-mmphi-sensitivity}

\end{refremark}

Our goal is to use the Morris-Mitchell criterion as an additional
objective in the multi-objective optimization problem to improve the
quality of the coverage of the input space. Therefore, it is crucial to
understand how $\Phi^{\ast}$ behaves when new points are added to an
existing design. First experimental observations motivated a rigorous
analysis of the sensitivity of $\Phi^{\ast}$ to added points, which is
presented in Section~\ref{sec-mm-monoticity-analysis}.

\subsection{Analysis of Monotonicity}\label{sec-mm-monoticity-analysis}

This section analyzes conditions under which $\Phi^{\ast}$ decreases
or increases when a new point is added to an existing design.

\begin{theorem}[Monotonicity behavior of
$\Phi^{\ast}$]\protect\hypertarget{thm-sensitivity-mmphi}{}\label{thm-sensitivity-mmphi}

The intensified Morris-Mitchell criterion $\Phi^{\ast}$ decreases if
and only if the average $d^{-q}$ contribution of the new point to the
existing points (weighted by the new distance multiplicities $J_k'$),
i.e.,

\[
\frac{1}{n}\sum_{k=1}^{m'} J_k' (d_k')^{-q},
\]

is strictly below the current design's overall average pair contribution
(weighted by the existing distance multiplicities $J_j$), i.e.,

\[
\frac{1}{\binom{n}{2}} \sum_{j=1}^{m} J_j d_j^{-q}.
\]

\end{theorem}

\begin{proof}
Let $X_p = \{x_1, \ldots, x_n\}$ be a design. According to
Definition~\ref{def-mmphi}, let $d_1 < d_2 < \dots < d_m$ be the
distinct distances between all pairs of points in $X_p$, and let
$J_j$ be the number of pairs separated by distance $d_j$. The
unnormalized $q$-th power of the default criterion
(Equation~\ref{eq-mmphi}) is: \[
\Phi^q = \Phi_q^q(X_p) = \sum_{j=1}^{m} J_j d_j^{-q}
\] where $\sum_{j=1}^m J_j = \binom{n}{2}$.

Let $x_{n+1} \notin X_p$ be a new point. The insertion of $x_{n+1}$
creates $n$ new point pairs. Let the distinct distances between
$x_{n+1}$ and the existing points in $X_p$ be
$d_1' < d_2' < \dots < d_{m'}'$, with corresponding multiplicities
$J_1', J_2', \dots, J_{m'}'$ such that $\sum_{k=1}^{m'} J_k' = n$.

We define the cross-term sum, which can be interpreted as the
\emph{interaction energy} of the new point with the existing design, as:
\[
\Delta(x_{n+1}, X_p) := \sum_{k=1}^{m'} J_k' (d_k')^{-q} > 0.
\]

The intensified Morris-Mitchell criterion $\Phi_q^{\ast}$ normalizes
the sum by the total number of pairs $M$. For the new design
$X_p \cup \{x_{n+1}\}$, the total number of pairs is
$\binom{n+1}{2}$. Thus, the intensive criterion evaluates to:

\begin{align*}
\Phi^{\ast q}(X_p \cup \{x_{n+1}\}) & = \frac{1}{\binom{n+1}{2}}\left[\sum_{j=1}^{m} J_j d_j^{-q} + \sum_{k=1}^{m'} J_k' (d_k')^{-q}\right] \\
& = \frac{2}{n(n+1)}\left[\Phi^q(X_p) + \Delta(x_{n+1}, X_p)\right].
\end{align*}

The condition for a strict decrease,
$\Phi^{\ast q}(X_p \cup \{x_{n+1}\}) < \Phi^{\ast q}(X_p)$, is
equivalent to:

\[
\frac{2\left[\Phi^q(X_p) + \Delta(x_{n+1}, X_p)\right]}{n(n+1)} < \frac{2\Phi^q(X_p)}{n(n-1)}.
\]

Dividing by $2/n > 0$ and rearranging gives:

\begin{align*}
(n-1)\left[\Phi_q^q(X_p) + \Delta(x_{n+1}, X_p)\right] & < (n+1)\Phi_q^q(X_p) \Leftrightarrow\\
(n-1)\Delta(x_{n+1}, X_p) < 2\Phi_q^q(X_p) & = n(n-1)\Phi_q^{\ast q}(X_p).
\end{align*}

Since $\Phi^q(X_p) = \binom{n}{2}\Phi^{\ast q}(X_p)$, we obtain the
necessary and sufficient criterion for a decrease:

\[
\frac{1}{n}\sum_{k=1}^{m'} J_k' (d_k')^{-q} < \frac{1}{\binom{n}{2}} \sum_{j=1}^{m} J_j d_j^{-q} = \Phi^{\ast q}(X_p).
\]
\end{proof}

\begin{example}[Adding an optimal point increases
$\Phi^{\ast}$]\protect\hypertarget{exm-counter}{}\label{exm-counter}

Let $S = [0,1]$, i.e., the one-dimensional unit interval,
$X_p = \{0,1\}$, and $q = 2$. There is $m=1$ distinct distance
$d_1 = 1$, with multiplicity $J_1 = 1$. \[
\Phi_2^{\ast 2}(X_p) = \frac{1}{\binom{2}{2}} \times J_1 d_1^{-2} = 1 \times 1 \times 1^{-2} = 1.
\]

We add $x_3 = 0.5$, which is geometrically optimal for a three-point
equidistant design. This can be computed in \texttt{spotoptim} as
follows:

\end{example}

\protect\phantomsection\label{mmphi-intensive-example-2}
\begin{Shaded}
\begin{Highlighting}[]
\NormalTok{X2 }\OperatorTok{=}\NormalTok{ np.array([}
\NormalTok{    [}\FloatTok{0.0}\NormalTok{],}
\NormalTok{    [}\FloatTok{1.0}\NormalTok{]}
\NormalTok{])}
\NormalTok{quality, J, d }\OperatorTok{=}\NormalTok{ mmphi\_intensive(X2, q}\OperatorTok{=}\DecValTok{2}\NormalTok{, p}\OperatorTok{=}\DecValTok{2}\NormalTok{)}
\BuiltInTok{print}\NormalTok{(}\SpecialStringTok{f"Quality (Phi\_q\_intensive): }\SpecialCharTok{\{}\NormalTok{quality}\SpecialCharTok{\}}\SpecialStringTok{"}\NormalTok{)}
\BuiltInTok{print}\NormalTok{(}\SpecialStringTok{f"Multiplicities (J): }\SpecialCharTok{\{}\NormalTok{J}\SpecialCharTok{\}}\SpecialStringTok{"}\NormalTok{)}
\BuiltInTok{print}\NormalTok{(}\SpecialStringTok{f"Distinct Distances (d): }\SpecialCharTok{\{}\NormalTok{d}\SpecialCharTok{\}}\SpecialStringTok{"}\NormalTok{)}
\NormalTok{new\_point }\OperatorTok{=}\NormalTok{ np.array([}\FloatTok{0.5}\NormalTok{])}
\NormalTok{new\_quality, new\_J, new\_d }\OperatorTok{=}\NormalTok{ mmphi\_intensive\_update(X2, new\_point, J, d, q}\OperatorTok{=}\DecValTok{2}\NormalTok{, p}\OperatorTok{=}\DecValTok{2}\NormalTok{)}
\BuiltInTok{print}\NormalTok{(}\SpecialStringTok{f"Updated Quality: }\SpecialCharTok{\{}\NormalTok{new\_quality}\SpecialCharTok{\}}\SpecialStringTok{"}\NormalTok{)}
\BuiltInTok{print}\NormalTok{(}\SpecialStringTok{f"Updated Multiplicities (J): }\SpecialCharTok{\{}\NormalTok{new\_J}\SpecialCharTok{\}}\SpecialStringTok{"}\NormalTok{)}
\BuiltInTok{print}\NormalTok{(}\SpecialStringTok{f"Updated Distinct Distances (d): }\SpecialCharTok{\{}\NormalTok{new\_d}\SpecialCharTok{\}}\SpecialStringTok{"}\NormalTok{)}
\end{Highlighting}
\end{Shaded}

\begin{verbatim}
Quality (Phi_q_intensive): 1.0
Multiplicities (J): [1]
Distinct Distances (d): [1.]
Updated Quality: 1.7320508075688772
Updated Multiplicities (J): [2 1]
Updated Distinct Distances (d): [0.5 1. ]
\end{verbatim}

\begin{figure}

\centering{

\pandocbounded{\includegraphics[keepaspectratio]{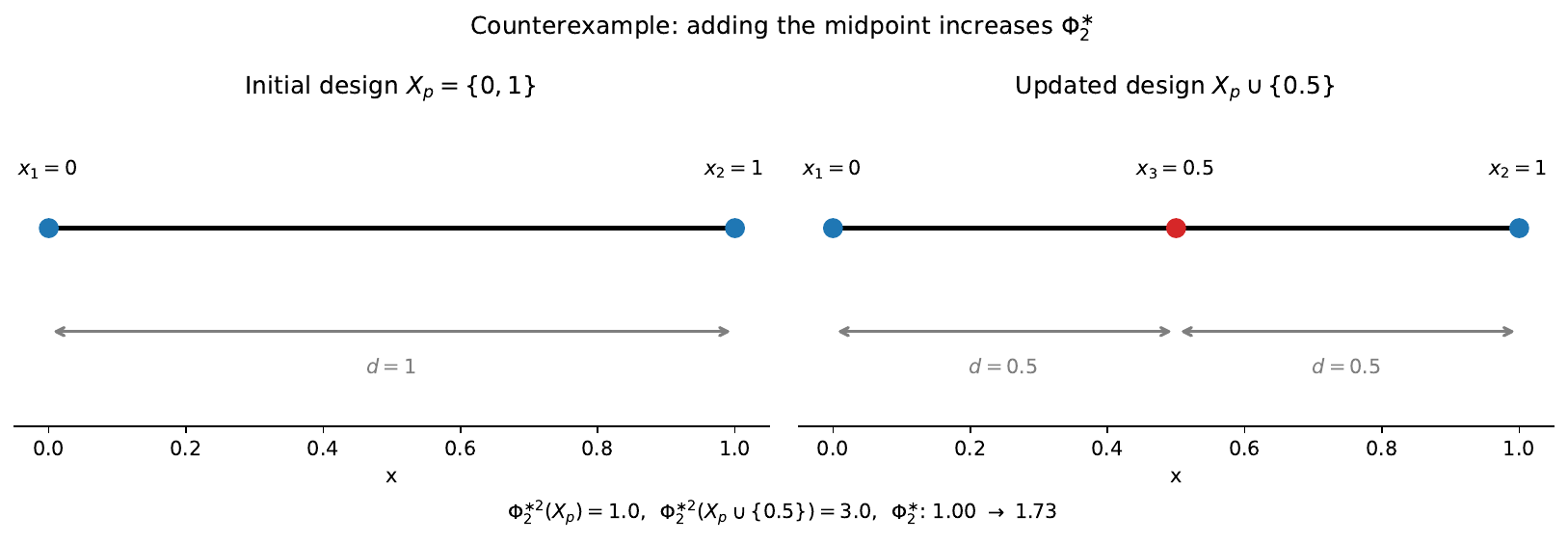}}

}

\caption{\label{fig-mmphi-intensive-update-example-counter}Counterexample:
adding the midpoint increases $\Phi_2^{\ast}$. Left panel: initial
design $P = \{0, 1\}$. Right panel: updated design
$P \cup \{0.5\}$.}

\end{figure}

The situation is illustrated in
Figure~\ref{fig-mmphi-intensive-update-example-counter}.

This introduces $m'=1$ new distinct distance $d_1' = 0.5$, which
occurs with multiplicity $J_1' = 2$ (since the distance to $0$ and
the distance to $1$ are both $0.5$). The cross-term is: \[
\Delta(0.5, P) = J_1' (d_1')^{-2} = 2 \times (0.5)^{-2} = 2 \times 4 = 8.
\]

The condition check \[
\frac{\Delta}{n} = \frac{8}{2} = 4 < \Phi_2^{\ast 2}(P) = 1
\] is false. Therefore,

\[
\Phi_2^{\ast 2}(P \cup \{0.5\}) = \frac{1}{3}(1+8) = 3 > 1 = \Phi_2^{\ast 2}(P).
\]

Hence $\Phi_2^{\ast}$ increases from $1$ to
$\sqrt{3} \approx 1.73$, although the new point is optimally placed.

\subsubsection{Connection to Potential
Theory}\label{connection-to-potential-theory}

\begin{theorem}[Scaling behavior of $\Phi$ and $\Phi^{\ast}$ for
optimal
designs]\protect\hypertarget{thm-mmphi-scaling}{}\label{thm-mmphi-scaling}

For an optimal $n$-point design on $[0,1]^k$ (e.g., a grid with
spacing $n^{-1/k}$), potential theory (classical Riesz energy theory)
(Björck 1956; Hardin and Saff 2005) yields the asymptotic growth of the
default Morris-Mitchell sum $\Phi(n)$ as $n \to \infty$:

\begin{equation}\protect\phantomsection\label{eq-mmphi-scaling}{
\Phi_q^q(n) = \sum_{j=1}^{m} J_j d_j^{-q} \sim  C_{q,k} \times n^{1+q/k},
\quad n \to \infty, \quad q > k.
}\end{equation}

With $M = \binom{n}{2} \sim n^2/2$, the intensive criterion
$\Phi_q^{\ast}$ asymptotically behaves as:

\[
\Phi_q^{\ast q}(n)
= \frac{\sum_{j=1}^{m} J_j d_j^{-q}}{\binom{n}{2}}
\sim \frac{ C_{q,k} \times n^{1+q/k}}{n^2/2}
= 2 \times C_{q,k} \times n^{q/k-1},
\]

where $C_{q,k}$ is a constant depending on $q$ and $k$.

\end{theorem}

\begin{proof}
Consider Equation 9 in Hardin and Saff (2005), the minimal Riesz
$s$-energy of $N$ points on a compact set $A \subset \mathbb{R}^d$
behaves asymptotically as:

\begin{equation}\protect\phantomsection\label{eq-riesz-asymptotic}{
\lim_{N \to \infty} \frac{E_s(A,N)}{N^{1+s/d}} = \frac{C_{s,d}}{H_d(A)^{s/d}},
}\end{equation}

where $E_s(A,N)$ is the minimal Riesz $s$-energy of $N$ points on
a compact set $A \subset \mathbb{R}^d$, $H_d(A)$ is the
$d$-dimensional Hausdorff measure of $A$, and $C_{s,d}$ is a
constant depending on $s$ and $d$. For the unit cube in
$\mathbb{R}^k$, we have $H_k([0,1]^k) = 1$.

Using Equation~\ref{eq-riesz-asymptotic}, the proof is basically a
substitution of the variables from the Riesz energy setting to the
Morris-Mitchell setting: The Riesz parameter $s$ corresponds to $q$,
the dimension $d$ corresponds to $k$. The number of points $N$
corresponds to $n$. The continuous $s$-energy (Riesz energy) of a
configuration

\[
E_s(A,N)=\sum_{i \neq j} |x_i - x_j|^{-s}
\]

sums over all \emph{ordered} pairs, so $E_s = 2\Phi_q^q$ (the
Morris-Mitchell sum counts \emph{unordered} pairs). Substituting the
translated variables into Equation~\ref{eq-riesz-asymptotic} yields
$E_q \sim C_{q,k} \times n^{1+q/k}$ and hence the asymptotic growth of
$\Phi_q^q(n)$ as $n \to \infty$ for optimal designs on $[0,1]^k$:
\[
\Phi_q^q(n) = \tfrac{1}{2}E_q([0,1]^k,n) \sim  \tfrac{1}{2}C_{q,k} \times n^{1+q/k},
\quad n \to \infty, \quad q > k.
\] Since the factor $\tfrac{1}{2}$ is a constant that can be absorbed
into $C_{q,k}$, we retain the scaling form
$\Phi_q^q(n) \sim C_{q,k}' \times n^{1+q/k}$ (writing $C_{q,k}'$ for
$\tfrac{1}{2}C_{q,k}$) as stated in Equation~\ref{eq-mmphi-scaling}.
\end{proof}

\begin{refremark}[The constant $C_{q,k}$]
Computing the exact value of $C_{q,k}$ is a notoriously difficult
problem in discrete geometry (Hardin et al. 2019) and beyond the scope
of this work. The exact value of $C_{q,k}$ is not required for the
scaling analysis, as it only affects the constant factor in the
asymptotic growth.

\label{rem-cqk}

\end{refremark}

\begin{refremark}[The case $q=k$]
If $q=k$ holds by coincidence, $\Phi_q^{\ast}(n)$ grows as
$\left(\log(n)\right)^{1/q}$ (Hardin et al. 2019).

\label{rem-qk}

\end{refremark}

In practice, one typically chooses $q \gg k$ (e.g., $q=10$,
$k=2$), so the intensive criterion $\Phi_q^{\ast}$ systematically
increases with $n$ even for optimal designs. Normalization by the
total pairs $M$ removes only the $n^2$ factor from the combinatorial
growth, but not the remaining $n^{q/k-1}$ geometric scaling factor.

\begin{example}[Scaling example: $\Phi_q^{\ast}$ increases with $n$
for optimal
designs]\protect\hypertarget{exm-scaling}{}\label{exm-scaling}

Verification example ($k=1$, $q=2$, equidistant points):

\begin{itemize}
\tightlist
\item
  $n=2$: $\Phi_2^{\ast} = 1.00$
\item
  $n=3$: $\Phi_2^{\ast} = \sqrt{3} \approx 1.73$
\item
  $n=4$: $\Phi_2^{\ast} = \sqrt{32.5/6} \approx 2.33$
\end{itemize}

This is strictly increasing, although each design is optimal for its
size.

\end{example}

A size-invariant normalization must eliminate the factor $n^{1+q/k}$
in Equation~\ref{eq-mmphi-scaling}. This leads to the corrected
Morris-Mitchell criterion.

\begin{definition}[Corrected Morris-Mitchell
Criterion]\protect\hypertarget{def-mmphi-corrected}{}\label{def-mmphi-corrected}

The corrected Morris-Mitchell criterion is defined as the original
criterion normalized by $n^{1+q/k}$. For a design $X$, it will be
denoted as $\hat{\Phi}_q$: \[
\hat{\Phi} =
\hat{\Phi}_q(X)
= \left(\frac{\sum_{j=1}^{m} J_j d_j^{-q}}{n^{1+q/k}}\right)^{1/q}
= \frac{\Phi_q(X)}{n^{1/q+1/k}}.
\]

\end{definition}

For optimal designs, $X_p$, the corrected criterion $\hat{\Phi}_q$
is asymptotically constant as $n \to \infty$:
$\hat{\Phi}_q(X_p) \to C_{q,k}^{1/q} = \text{const}$, and therefore \[
\hat{\Phi}_q(X_p \cup \{x_{n+1}^{\ast}\})
\lesssim
\hat{\Phi}_q(X_p)
\] for sufficiently large $n$ and optimally placed $x_{n+1}^{\ast}$.
This yields the Lemma~\ref{lem-mmphi-scaling}:

\begin{lemma}[Lemma: Scaling behavior of Morris-Mitchell
criteria]\protect\hypertarget{lem-mmphi-scaling}{}\label{lem-mmphi-scaling}

$\Phi_q^{\ast}$ does not possess the desired size-invariance property
in full generality. Normalization by $M = \binom{n}{2}$ is only
asymptotically correct when $q = k$ holds by coincidence. For
$q > k$, a correction by $n^{1+q/k}$ is mathematically required ---
at the cost of explicit dependence on the dimension $k$.

\end{lemma}

$\hat{\Phi}_q$ is asymptotically size-invariant, but not monotonically
decreasing when optimally placed points are added for all finite $n$.
The desired property --- monotone decrease when optimally placed points
are added --- cannot be achieved by any normalization of the form
$\Phi_q^q / f(n)$ for all finite $n$. This is because in a bounded
space the minimum distance necessarily decreases when going from $n$
to $n+1$ points, and the resulting small distances are heavily
penalized by the $d^{-q}$ term. The issue lies in the underlying
concept: as a maximin criterion, $\Phi_q$ measures packing quality
(how far points are apart). In a bounded space, the minimum distance
must decrease as $n$ increases. By contrast, the dual concept ---
covering quality (how well the space is covered using a minimax
criterion) --- has natural monotonicity.

\subsubsection{\texorpdfstring{\texttt{mmphi\_corrected}: the corrected
Morris-Mitchell criterion in
\texttt{spotoptim}}{mmphi\_corrected: the corrected Morris-Mitchell criterion in spotoptim}}\label{mmphi_corrected-the-corrected-morris-mitchell-criterion-in-spotoptim}

\begin{example}[Usage of
\texttt{mmphi\_corrected}]\protect\hypertarget{exm-mmphi-corrected}{}\label{exm-mmphi-corrected}

We create a simple 3-point sampling plan in 2D and calculate the
corrected Morris-Mitchell criterion with q=2, using Euclidean distances
(p=2).

\end{example}

\protect\phantomsection\label{mmphi-corrected-example}
\begin{Shaded}
\begin{Highlighting}[]
\NormalTok{X3 }\OperatorTok{=}\NormalTok{ np.array([}
\NormalTok{    [}\FloatTok{0.0}\NormalTok{, }\FloatTok{0.0}\NormalTok{],}
\NormalTok{    [}\FloatTok{0.5}\NormalTok{, }\FloatTok{0.5}\NormalTok{],}
\NormalTok{    [}\FloatTok{1.0}\NormalTok{, }\FloatTok{1.0}\NormalTok{]}
\NormalTok{])}
\NormalTok{quality, J, d }\OperatorTok{=}\NormalTok{ mmphi\_corrected(X3, q}\OperatorTok{=}\DecValTok{2}\NormalTok{, p}\OperatorTok{=}\DecValTok{2}\NormalTok{)}
\BuiltInTok{print}\NormalTok{(}\SpecialStringTok{f"Quality (Phi\_q\_corrected): }\SpecialCharTok{\{}\NormalTok{quality}\SpecialCharTok{\}}\SpecialStringTok{"}\NormalTok{)}
\BuiltInTok{print}\NormalTok{(}\SpecialStringTok{f"Multiplicities (J): }\SpecialCharTok{\{}\NormalTok{J}\SpecialCharTok{\}}\SpecialStringTok{"}\NormalTok{)}
\BuiltInTok{print}\NormalTok{(}\SpecialStringTok{f"Distinct Distances (d): }\SpecialCharTok{\{}\NormalTok{d}\SpecialCharTok{\}}\SpecialStringTok{"}\NormalTok{)}
\end{Highlighting}
\end{Shaded}

\begin{verbatim}
Quality (Phi_q_corrected): 0.7071067811865475
Multiplicities (J): [2 1]
Distinct Distances (d): [0.70710678 1.41421356]
\end{verbatim}

The returned \texttt{quality} is the $\hat{\Phi}$ value. The array
\texttt{J} contains the multiplicities for each distinct distance, and
\texttt{d} contains the distinct distances found in the design. There
are three points, resulting in three distances, $d_{12}$, $d_{13}$,
and $d_{23}$. The smallest distance is $d_{12} = d_{23}$, and it
occurs twice, so $J_1 = 2$. The largest distance is $d_{13}$, and it
occurs once, so $J_2 = 1$.

\begin{example}[Usage of
\texttt{mmphi\_corrected\_update}]\protect\hypertarget{exm-mmphi-corrected-update}{}\label{exm-mmphi-corrected-update}

Assume we have the state from the previous example: X, quality, J, d,
and we want to add a new point.

\end{example}

\protect\phantomsection\label{mmphi-corrected-update-example-onepoint}
\begin{Shaded}
\begin{Highlighting}[]
\NormalTok{new\_point }\OperatorTok{=}\NormalTok{ np.array([}\FloatTok{0.1}\NormalTok{, }\FloatTok{0.1}\NormalTok{])}
\NormalTok{new\_quality, new\_J, new\_d }\OperatorTok{=}\NormalTok{ mmphi\_corrected\_update(X3, new\_point, J, d, q}\OperatorTok{=}\DecValTok{2}\NormalTok{, p}\OperatorTok{=}\DecValTok{2}\NormalTok{)}
\BuiltInTok{print}\NormalTok{(}\SpecialStringTok{f"Updated Quality: }\SpecialCharTok{\{}\NormalTok{new\_quality}\SpecialCharTok{\}}\SpecialStringTok{"}\NormalTok{)}
\end{Highlighting}
\end{Shaded}

\begin{verbatim}
Updated Quality: 1.907915812323379
\end{verbatim}

The function \texttt{mmphi\_corrected\_update} returns the updated
criterion, multiplicities, and distinct distances, which can be used for
the next update step.

\begin{example}[Adding an optimal point increases
$\hat{\Phi}$]\protect\hypertarget{exm-counter-2}{}\label{exm-counter-2}

We use the same setup as in Example~\ref{exm-counter}, but now we use
the corrected Morris-Mitchell criterion.

Let $S = [0,1]$, i.e., the one-dimensional ($k=1$) unit interval,
$P = \{0,1\}$, and $q = 2$. There is $m=1$ distinct distance
$d_1 = 1$, with multiplicity $J_1 = 1$. The normalization is
computed as \$ n\^{}\{1.0 + q / k\} = 2\^{}\{1.0 + 2 / 1\} = 2\^{}\{3\}
= 8\$ \[
\hat{\Phi}_2^{2}(P) = \frac{1}{n^{1.0 + q / k}} \times J_1 d_1^{-2} = \frac{1}{8} \times 1 \times 1^{-2} = 0.125 \Longrightarrow \hat{\Phi}_2(P) = \sqrt{0.125} \approx 0.35355.
\]

Add $x_3 = 0.5$, which is geometrically optimal for a 3-point
equidistant design. This can be computed in \texttt{spotoptim} as
follows:

\end{example}

\protect\phantomsection\label{mmphi-corrected-example-2}
\begin{Shaded}
\begin{Highlighting}[]
\NormalTok{X2 }\OperatorTok{=}\NormalTok{ np.array([}
\NormalTok{    [}\FloatTok{0.0}\NormalTok{],}
\NormalTok{    [}\FloatTok{1.0}\NormalTok{]}
\NormalTok{])}
\NormalTok{quality, J, d }\OperatorTok{=}\NormalTok{ mmphi\_corrected(X2, q}\OperatorTok{=}\DecValTok{2}\NormalTok{, p}\OperatorTok{=}\DecValTok{2}\NormalTok{)}
\BuiltInTok{print}\NormalTok{(}\SpecialStringTok{f"Quality (Phi\_q\_corrected): }\SpecialCharTok{\{}\NormalTok{quality}\SpecialCharTok{\}}\SpecialStringTok{"}\NormalTok{)}
\BuiltInTok{print}\NormalTok{(}\SpecialStringTok{f"Multiplicities (J): }\SpecialCharTok{\{}\NormalTok{J}\SpecialCharTok{\}}\SpecialStringTok{"}\NormalTok{)}
\BuiltInTok{print}\NormalTok{(}\SpecialStringTok{f"Distinct Distances (d): }\SpecialCharTok{\{}\NormalTok{d}\SpecialCharTok{\}}\SpecialStringTok{"}\NormalTok{)}
\NormalTok{new\_point }\OperatorTok{=}\NormalTok{ np.array([}\FloatTok{0.5}\NormalTok{])}
\NormalTok{new\_quality, new\_J, new\_d }\OperatorTok{=}\NormalTok{ mmphi\_corrected\_update(X2, new\_point, J, d, q}\OperatorTok{=}\DecValTok{2}\NormalTok{, p}\OperatorTok{=}\DecValTok{2}\NormalTok{)}
\BuiltInTok{print}\NormalTok{(}\SpecialStringTok{f"Updated Quality: }\SpecialCharTok{\{}\NormalTok{new\_quality}\SpecialCharTok{\}}\SpecialStringTok{"}\NormalTok{)}
\BuiltInTok{print}\NormalTok{(}\SpecialStringTok{f"Updated Multiplicities (J): }\SpecialCharTok{\{}\NormalTok{new\_J}\SpecialCharTok{\}}\SpecialStringTok{"}\NormalTok{)}
\BuiltInTok{print}\NormalTok{(}\SpecialStringTok{f"Updated Distinct Distances (d): }\SpecialCharTok{\{}\NormalTok{new\_d}\SpecialCharTok{\}}\SpecialStringTok{"}\NormalTok{)}
\end{Highlighting}
\end{Shaded}

\begin{verbatim}
Quality (Phi_q_corrected): 0.3535533905932738
Multiplicities (J): [1]
Distinct Distances (d): [1.]
Updated Quality: 0.5773502691896257
Updated Multiplicities (J): [2 1]
Updated Distinct Distances (d): [0.5 1. ]
\end{verbatim}

\begin{figure}

\centering{

\pandocbounded{\includegraphics[keepaspectratio]{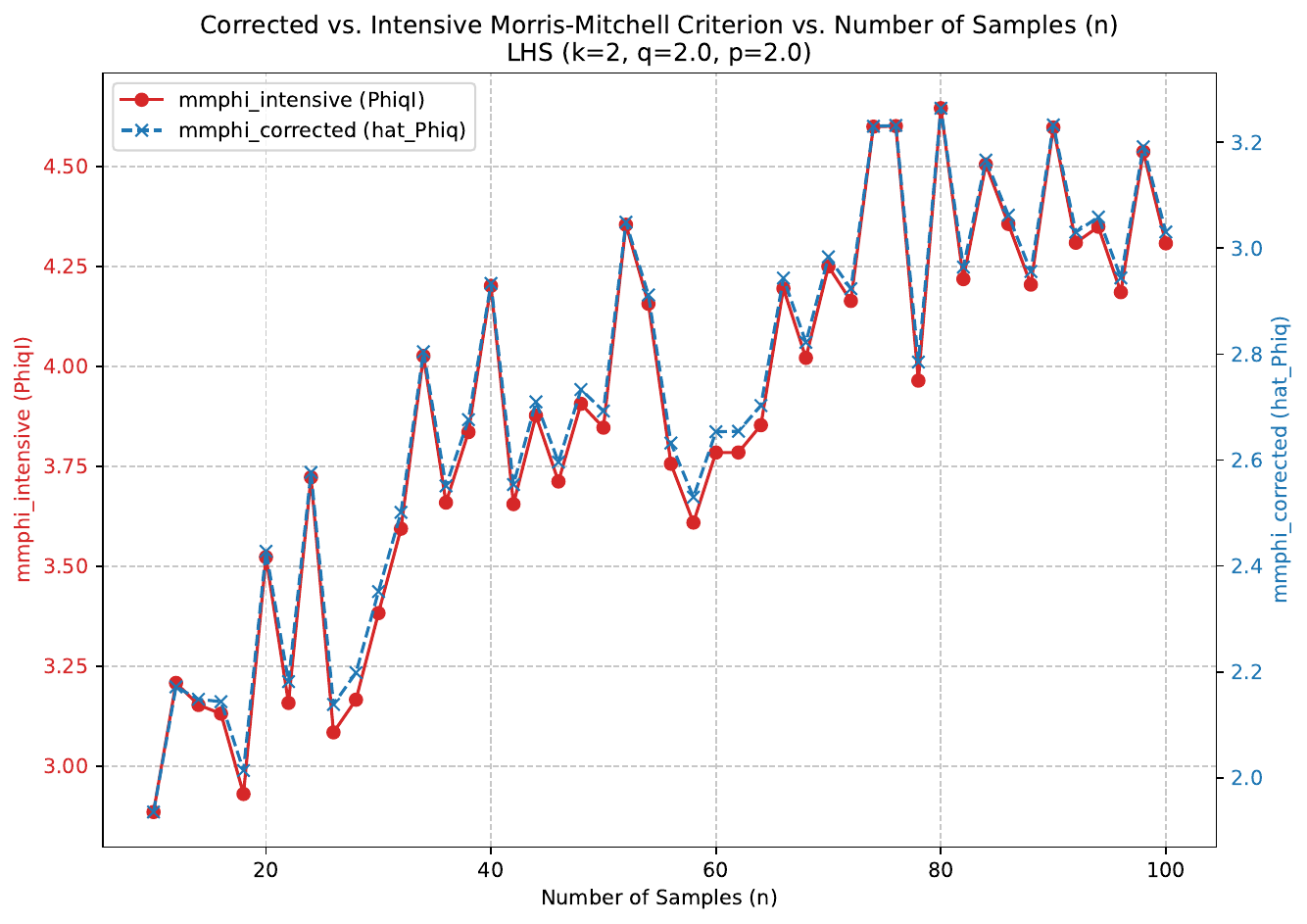}}

}

\caption{\label{fig-mm-corrected-vs-n-lhs}Morris-Mitchell criterion and
corrected criterion vs.~number of samples (n) for LHS designs. Note the
different scales of the y-axes.}

\end{figure}

Figure~\ref{fig-mm-corrected-vs-n-lhs} shows the results of the
intensified Morris-Mitchell criterion and its corrected version versus
the number of samples ($n$) for LHS designs in a 2-dimensional space.
It was generated with the function
\texttt{plot\_mmphi\_corrected\_vs\_n\_lhs} from the \texttt{spotoptim}
package. Figure~\ref{fig-mm-corrected-vs-n-lhs} shows that both criteria
show a similar behavior, because they only differ by the normalization
factor.

\begin{figure}
\centering{
\pandocbounded{\includegraphics[width=0.8\linewidth,keepaspectratio]{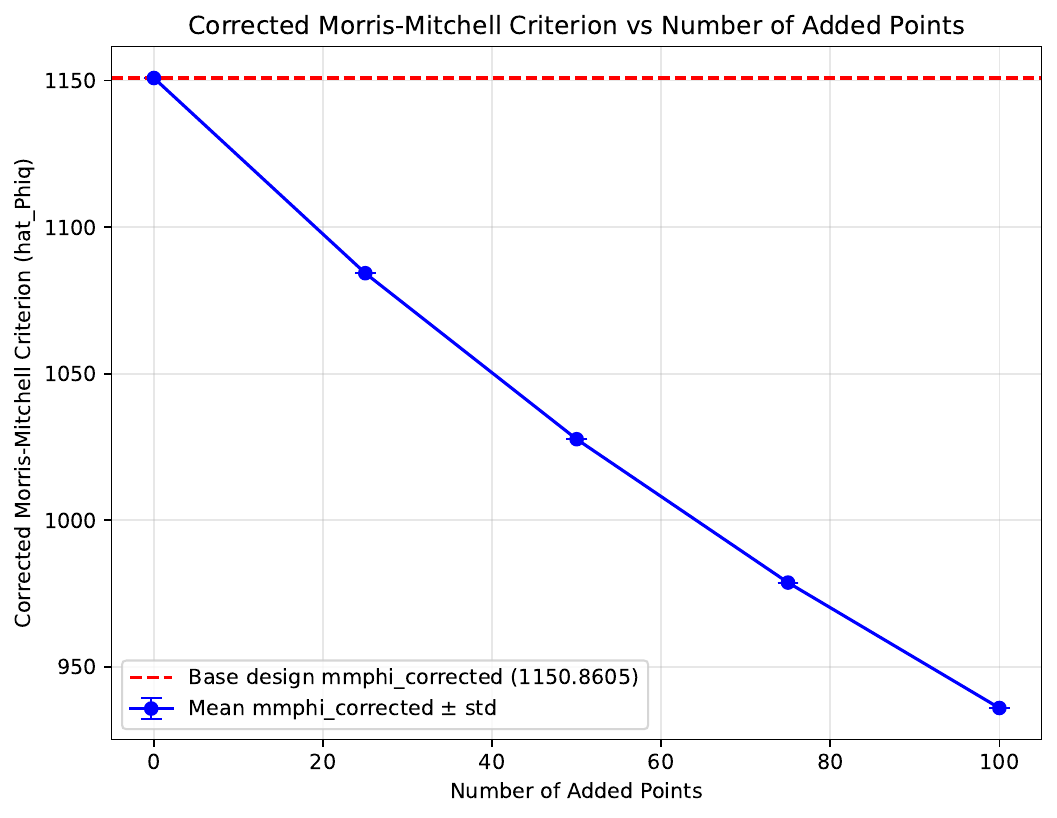}}
}
\caption{\label{fig-mmphi-corrected-plot-x-call}Compressor data. The
corrected Morris-Mitchell criterion vs.~number of added points.
Randomized points are added to the entire dataset, which results in a
decrease of the criterion.}
\end{figure}

Figure~\ref{fig-mmphi-corrected-plot-x-call} shows the Morris-Mitchell
criterion versus the number of added points. Randomized points are added
to the entire dataset.

\begin{figure}
\centering{
\pandocbounded{\includegraphics[width=0.8\linewidth,keepaspectratio]{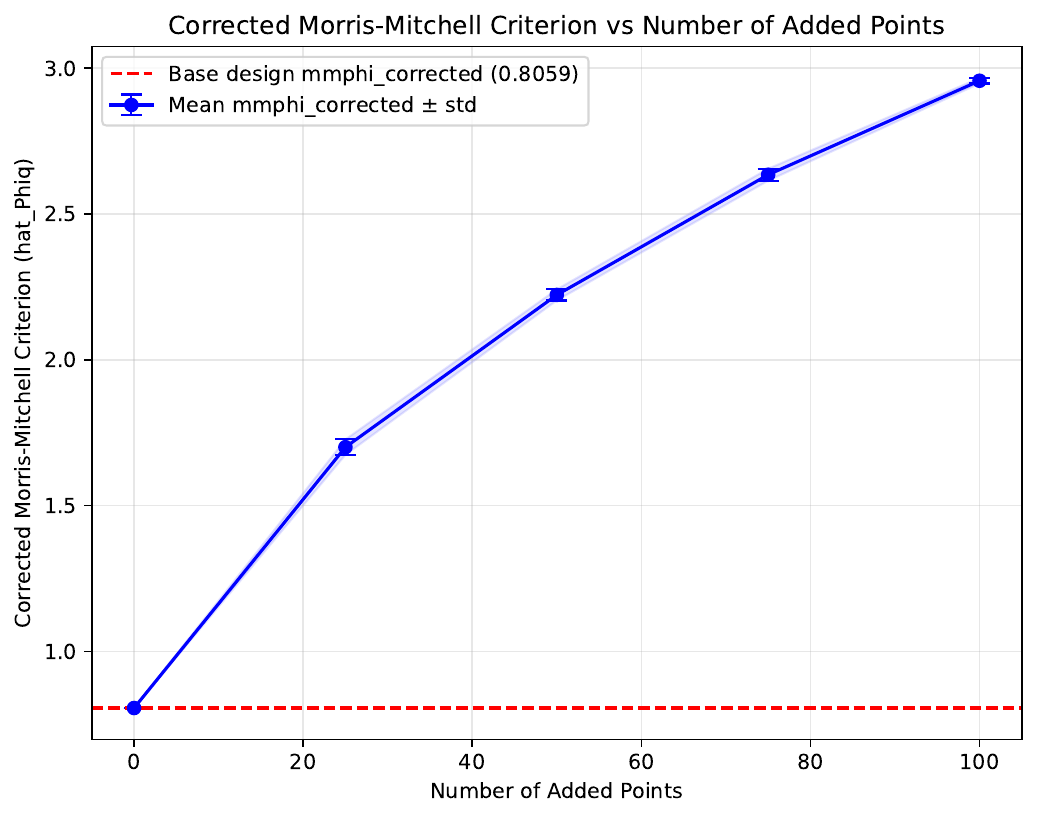}}
}
\caption{\label{fig-mmphi-corrected-plot-x-search-scaled1}Optimal
design. The corrected Morris-Mitchell criterion vs.~number of added
points. Randomized points are added to an LHD, which results in an
increase of the criterion.}

\end{figure}

In contrast to the previous example
(Figure~\ref{fig-mmphi-corrected-plot-x-call}), the corrected
Morris-Mitchell criterion is calculated for a planned and optimized
experimental design, i.e., random points are added to an existing design
$X_p$, which is nearly optimal.
Figure~\ref{fig-mmphi-corrected-plot-x-search-scaled1} shows the values
of the corrected Morris-Mitchell criterion $\hat{\Phi}_q$ after adding
random points to the existing design $X_p$.


To conclude our theoretical considerations, we present the relationship
between the intensified and corrected Morris-Mitchell criteria, stated
as Theorem~\ref{thm-mm-relationship}.

\begin{theorem}[Relationship between Intensified and Corrected
Morris-Mitchell
Criteria]\protect\hypertarget{thm-mm-relationship}{}\label{thm-mm-relationship}

Let $X = \{x_1, \ldots, x_n\}$ be a design of $n \ge 2$ points in a
$k$-dimensional space. Let $q \in \mathbb{Z}^+$ be the distance
weighting parameter (Riesz parameter), $d_j$ be the distinct pairwise
distances, and $J_j$ be the corresponding distance multiplicities.

Define the intensified Morris-Mitchell criterion, normalized by the
total number of pairs $\binom{n}{2}$, as:
\[ \Phi_q^{\ast}(X) = \left( \frac{2}{n(n-1)} \sum_{j=1}^m J_j d_j^{-q} \right)^{1/q} \]

Define the corrected Morris-Mitchell criterion, normalized by the
asymptotic Riesz energy scaling factor $n^{1+q/k}$ (derived from the
Poppy-seed bagel theorem (Wikipedia contributors 2025)), as:
\[ \hat{\Phi}_q(X) = \left( \frac{1}{n^{1+q/k}} \sum_{j=1}^m J_j d_j^{-q} \right)^{1/q} \]

Then, the relationship between $\hat{\Phi}_q(X)$ and
$\Phi_q^{\ast}(X)$ depends strictly on the relationship between the
parameter $q$ and the spatial dimension $k$:

\begin{itemize}
\tightlist
\item
  $q \ge k$ (hypersingular and critical regimes): For all $n \ge 2$,
  the corrected criterion is strictly smaller than the intensified
  criterion: \[ \hat{\Phi}_q(X) < \Phi_q^{\ast}(X) \]
\item
  $q < k$ (sub-harmonic regime): There exists a threshold size
  $N \in \mathbb{N}$ such that for all designs with $n > N$ points,
  the corrected criterion is strictly larger than the intensified
  criterion: \[ \hat{\Phi}_q(X) > \Phi_q^{\ast}(X) \]
\end{itemize}

\end{theorem}

\begin{proof}
To compare the two criteria, we evaluate the ratio of their $q$-th
powers, denoted as $R(n)$:
\[ R(n) = \frac{\hat{\Phi}_q^q(X)}{\Phi_q^{\ast q}(X)} = \frac{\left( \frac{1}{n^{1+q/k}} \sum_{j=1}^m J_j d_j^{-q} \right)}{\left( \frac{2}{n(n-1)} \sum_{j=1}^m J_j d_j^{-q} \right)} \]

Canceling the common unnormalized Morris-Mitchell sum
$\sum J_j d_j^{-q}$ yields:
\[ R(n) = \frac{n(n-1)}{2n^{1+q/k}} = \frac{n-1}{2n^{q/k}} \]

Since the function $f(x) = x^{1/q}$ is strictly monotonically
increasing for $x>0$, it holds that
$\hat{\Phi}_q(X) < \Phi_q^{\ast}(X)$ if and only if $R(n) < 1$, and
$\hat{\Phi}_q(X) > \Phi_q^{\ast}(X)$ if and only if $R(n) > 1$.

\begin{itemize}
\item
  Case 1: $q \ge k$: If $q \ge k$, the exponent ratio
  $\frac{q}{k} \ge 1$. Because $n \ge 2$, we have
  $n^{q/k} \ge n^1 = n$. Therefore, substituting this into the
  denominator of $R(n)$:
  \[ R(n) = \frac{n-1}{2n^{q/k}} \le \frac{n-1}{2n} \] Since
  $n - 1 < n$, it trivially follows that:
  \[ \frac{n-1}{2n} < \frac{n}{2n} = \frac{1}{2} < 1 \] Thus,
  $R(n) < 1$ is universally true for all $n \ge 2$. Consequently,
  $\hat{\Phi}_q(X) < \Phi_q^{\ast}(X)$.
\item
  Case 2: $q < k$: If $q < k$, the exponent ratio
  $\frac{q}{k} < 1$. Let $\frac{q}{k} = 1 - \epsilon$ for some
  $\epsilon > 0$. Rewriting $R(n)$ in terms of $\epsilon$:
  \[ R(n) = \frac{n-1}{2n^{1-\epsilon}} = \frac{n^\epsilon - n^{\epsilon-1}}{2} \]
  We analyze the limit of $R(n)$ as $n \to \infty$. Because
  $\epsilon > 0$, the term $n^\epsilon \to \infty$. Because
  $\epsilon - 1 < 0$, the term $n^{\epsilon-1} \to 0$.
  \[ \lim_{n \to \infty} R(n) = \lim_{n \to \infty} \frac{n^\epsilon - 0}{2} = \infty \]
  Since the sequence diverges to infinity, there must exist some
  $N \in \mathbb{N}$ such that for all $n > N$, $R(n) > 1$.
  Consequently, for sufficiently large designs,
  $\hat{\Phi}_q(X) > \Phi_q^{\ast}(X)$.
\end{itemize}

\end{proof}

\section{Desirability Functions}\label{sec-desirability}

The ``National Institute of Standards and Technology/SEMATECH e-Handbook
of Statistical Methods'' describes desirability functions as follows:

Desirability functions are among the most commonly used methods in
industry for optimizing processes with multiple objectives. This
approach is based on the idea that the ``quality'' of a product or
process with multiple quality characteristics is unacceptable if one of
these characteristics falls outside certain ``desired'' limits. The
method determines operating conditions that provide the ``most
desirable'' target values (National Institute of Standards and
Technology 2021).

The desirability function approach for simultaneous optimization of
multiple equations was originally proposed by Harington (1965).
Essentially, the approach consists of transforming the functions to a
common scale ($[0, 1]$), combining them using the geometric mean, and
optimizing the overall criterion. The equations can represent model
predictions or other equations (Kuhn 2016), (Bartz-Beielstein 2025b).
Desirability functions are popular, for example, in response surface
methodology Myers et al. (2016) as a method for simultaneously
optimizing a series of quadratic models. A response surface experiment
can use measurements on a series of outcomes. Instead of optimizing each
outcome separately, settings for the predictor variables are sought to
satisfy all outcomes simultaneously. Also in drug discovery, prediction
models can be created to relate the molecular structures of compounds to
properties of interest (e.g., absorption properties, efficacy, and
selectivity for the intended target). Given a set of prediction models
built using existing compounds, predictions can be made about a large
set of virtual compounds that have been designed but not necessarily
synthesized. Using the model predictions, a virtual compound can be
evaluated based on how well the model results match the required
properties. In this case, ranking compounds based on multiple endpoints
may be sufficient to meet the scientist's needs.

Originally, Harrington used exponential functions to quantify
desirability. We will use the simple discontinuous functions of
Derringer and Suich (1980). Suppose there are $R$ equations or
functions to be optimized simultaneously, denoted as $f_r(\vec{x})$
($r=1\ldots R$). For each of the $R$ functions, an individual
``desirability'' function is created that is high when $f_r(\vec{x})$
is at the desired level (e.g., a maximum, minimum, or target), and low
when $f_r(\vec{x})$ has an undesirable value. Derringer and Suich
(1980) proposed three forms of these functions corresponding to the type
of optimization objective. For maximizing $f_r(\vec{x})$, the function

\begin{equation}\protect\phantomsection\label{eq-Dmax}{
d_r^{\max}=
\begin{cases}
    0       &\text{if $f_r(\vec{x})< A$}    \\
    \left(\frac{f_r(\vec{x})-A}{B-A}\right)^{s}
            &\text{if $A \leq f_r(\vec{x}) \leq B$} \\
    1       &\text{if $f_r(\vec{x})>B$}
\end{cases}
}\end{equation}

can be used, where $A$, $B$, and $s$ are chosen by the
investigator. If the equation is to be minimized, they proposed the
function

\begin{equation}\protect\phantomsection\label{eq-Dmin}{d_r^{\min}=
\begin{cases}
    0       &\text{if $f_r(\vec{x})> B$}    \\
    \left(\frac{f_r(\vec{x})-B}{A-B}\right)^{s}
            &\text{if $A \leq f_r(\vec{x}) \leq B$} \\
    1       &\text{if $f_r(\vec{x})< A$}
\end{cases},
}\end{equation}

and for ``target-is-best'' situations (where the goal is to achieve a
specific target value $t_0$), they proposed the function

\begin{equation}\protect\phantomsection\label{eq-Dtarget}{
d_r^{\text{target}}=
\begin{cases}
    \left(\frac{f_r(\vec{x})-A}{t_0-A}\right)^{s_1}
            &\text{if $A \leq f_r(\vec{x}) \leq t_0$}   \\
    \left(\frac{f_r(\vec{x})-B}{t_0-B}\right)^{s_2}
            &\text{if $t_0 \leq f_r(\vec{x}) \leq B$}   \\
    0       &\text{otherwise}   \\
\end{cases}
}\end{equation}

These functions are on the same scale and are discontinuous at the
points $A$, $B$, and $t_0$. The values of $s$, $s_1$, or
$s_2$ can be chosen so that the desirability criterion is easier or
more difficult to satisfy. For example, if $s$ in
Equation~\ref{eq-Dmin} is chosen to be less than 1, $d_r^{\min}$ is
close to 1 even if the model $f_r(\vec{x})$ is not low. As the values
of $s$ move closer to 0, the desirability reflected by
Equation~\ref{eq-Dmin} becomes higher. Similarly, values of $s$
greater than 1 will make it more difficult to satisfy $d_r^{\min}$ in
terms of desirability. These scaling factors are useful when one
equation is more important than the others. It should be noted that any
function can be used to reflect the desirability of a model.

We are using the implementation of Bartz-Beielstein (2025a), which is
based on the original implementation in R of Kuhn (2016). The
\texttt{desirability} package (Kuhn 2016), which is written in the
statistical programming language \texttt{R}, contains \texttt{S3}
classes for multivariate optimization using the desirability function
approach of Harington (1965) with functional forms described by
Derringer and Suich (1980). It is available on CRAN, see
\url{https://cran.r-project.org/package=desirability}. A newer version,
the \texttt{desirability2} package, improves on the original
desirability package by enabling in-line computations that can be used
with dplyr pipelines (Kuhn 2025). It is also available on CRAN, see
\url{https://cran.r-project.org/web/packages/desirability2/index.html}.

\section{Surrogate-Model Based Optimization}\label{sec-optimization}

Because only a limited number of data points are available in practice,
a surrogate model is used to simulate the ground truth. This allows to
evaluate the objective function at any point in the design space and to
generate as many synthetic data points as desired. As discussed in
Section~\ref{sec-surrogate-models}, the random forest model showed the
best performance. Therefore, the random forest model, fitted on the
entire dataset is used as the ground truth. Using a fit on the entire
dataset is feasible here, because we want to use the complete
information available and we are not testing the performance of the
surrogate model.

The objective function evaluation during the optimization is performed
using a surrogate model trained on the entire dataset (\texttt{X\_full}
and \texttt{y\_full}). Since we are considering two-objectives,
\texttt{z8} and \texttt{z1}, two separate random forest models are
trained, one for each target variable: \texttt{z8\_ground\_rf} and
\texttt{z1\_ground\_rf}\footnote{As mentioned above, alternatively one random forest model could be trained for both target variables.}.

\section{\texorpdfstring{The Multi-Objective Case: Combining the two
Objectives \texttt{z8} and \texttt{z1} without
Morris-Mitchell}{The Multi-Objective Case: Combining the two Objectives z8 and z1 without Morris-Mitchell}}\label{sec-opt-wo-mm}

As the first objective function maximizing the first target variable,
\texttt{z8}, is used. The second objective is given by maximizing the
\texttt{z1} target variable. In general, if $x$ denotes a
$k$-dimensional input vector and $f_i(x)$ is the objective function
for the $i$-th target variable, the multi-objective optimization
problem is to simultaneously maximize all $p$ objectives: \[
\max_{x \in \mathcal{X}} \bigl(f_1(x), \ldots, f_p(x)\bigr).
\]

If $p=2$, we are facing a bi-objective problem. Using a desirability
function, the $p$ criteria are combined to obtain an overall
desirability \texttt{z}. The optimizer uses the desirability \texttt{z}
as a combination of the $p$ objective functions for the search:
\begin{equation}\protect\phantomsection\label{eq-overall-desirability}{
z = D_{\text{overall}}(f_1(x), \ldots, f_p(x)).
}\end{equation}

\subsection{Determination of Desirabilities for Target
Variables}\label{determination-of-desirabilities-for-target-variables}

For the desirability ranges, i.e., the intervals in which the
desirabilities are greater than zero, the boundaries \texttt{z\_min} and
\texttt{z\_max} are used. We first determine the boundaries for the
target variable \texttt{z8}, which is to be maximized. There are several
ways of setting the boundaries for the desirability functions. For
example, we can set the lower boundary to 90\% of the maximum value and
the upper boundary to 110\% of the maximum value. In our case, where the
data is normalized, this results in \texttt{z8\_min\ =\ 0.9} and
\texttt{z8\_max\ =\ 1.1}.

We do not choose this approach, because it results in a large, flat
(plateau-like) desirability function, which makes optimization
difficult. Instead, we use the shape parameters of the desirability
functions to adjust the steepness of the desirability functions.
Therefore, we can use the actual minimum and maximum values of the
target variables to set the boundaries for the desirability functions.
The maximum is extended by 10\% to allow for better
optimization\footnote{To implement this, we use the scale parameters
  \texttt{z8\_min\_multiplier=1.0} and \texttt{z8\_max\_multiplier=1.1}
  as well as \texttt{z1\_min\_multiplier} and
  \texttt{z1\_max\_multiplier=1.1}.}. Summarizing, we determine the
boundaries as follows: \texttt{low\ =\ z\_min},
\texttt{high\ =\ 1.1\ *\ z\_max}. Furthermore, we will use the scale
$s = 5$ for the desirability functions.

\begin{figure}
\centering{
\pandocbounded{\includegraphics[width=0.8\linewidth,keepaspectratio]{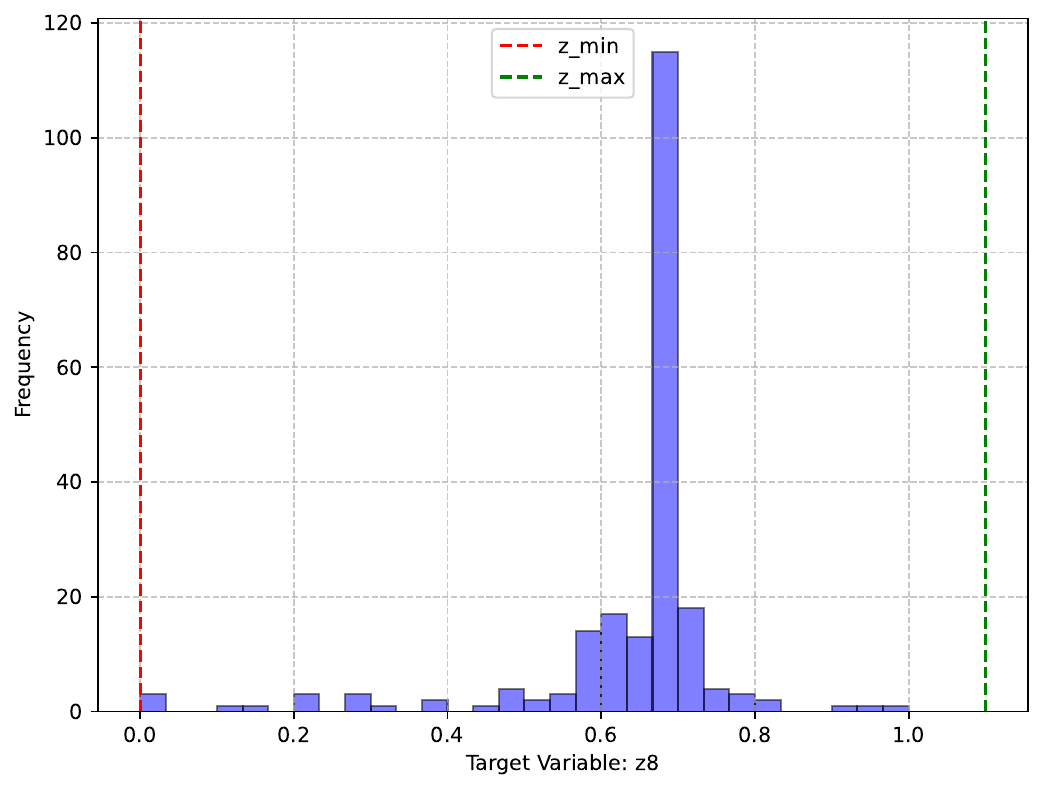}}
}
\caption{\label{fig-z-bounds-1dim-z8}Histogram of target variable with
desirability bounds. The red dashed line indicates the lower boundary
for the desirability function, and the green dashed line indicates the
upper boundary for the desirability function.}
\end{figure}

\begin{figure}
\centering{
\pandocbounded{\includegraphics[width=0.8\linewidth,keepaspectratio]{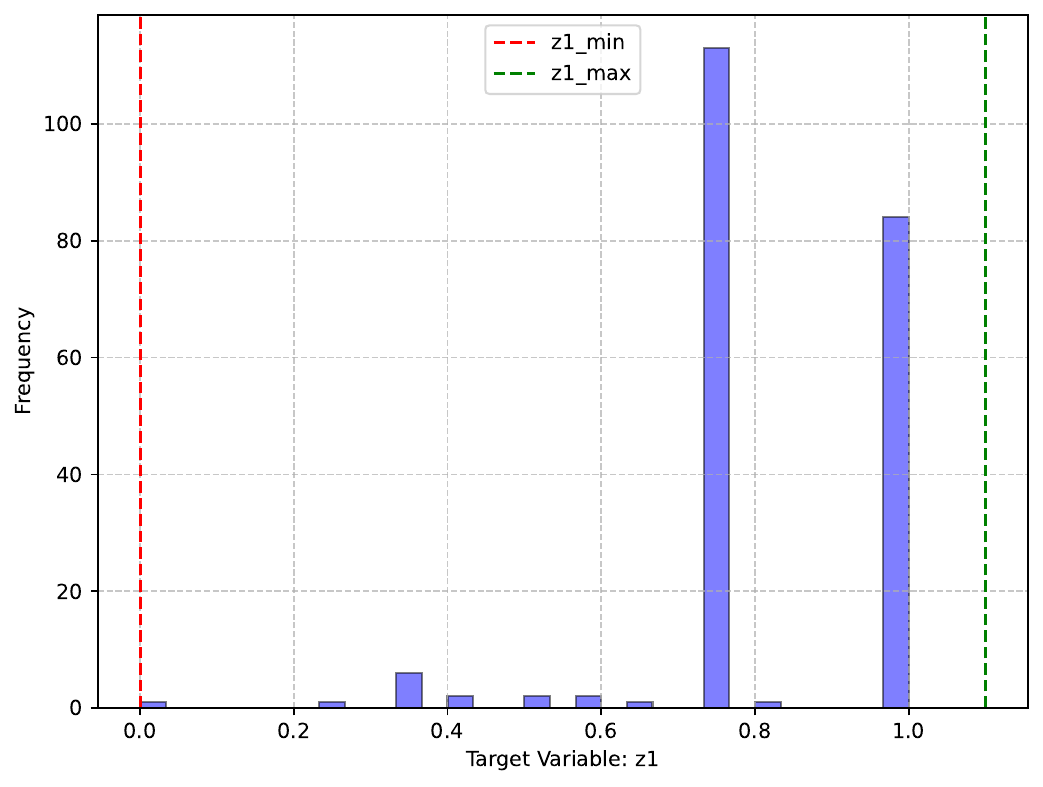}}
}
\caption{\label{fig-z-bounds-1dim-z1}Histogram of target variable with
desirability bounds. The red dashed line indicates the lower boundary
for the desirability function, and the green dashed line indicates the
upper boundary for the desirability function.}

\end{figure}

Figure~\ref{fig-z-bounds-1dim-z8} shows the histogram of the target
variable \texttt{z8} with the boundaries for the desirability indicated
as dashed lines. The boundaries for the target variable \texttt{z1},
which is also to be maximized, are set to identical values as for
\texttt{z8}. Figure~\ref{fig-z-bounds-1dim-z1} shows the histogram of
the target variable \texttt{z1} with the boundaries for the desirability
indicated as dashed lines.

\begin{figure}
\centering{
\pandocbounded{\includegraphics[width=0.8\linewidth,keepaspectratio]{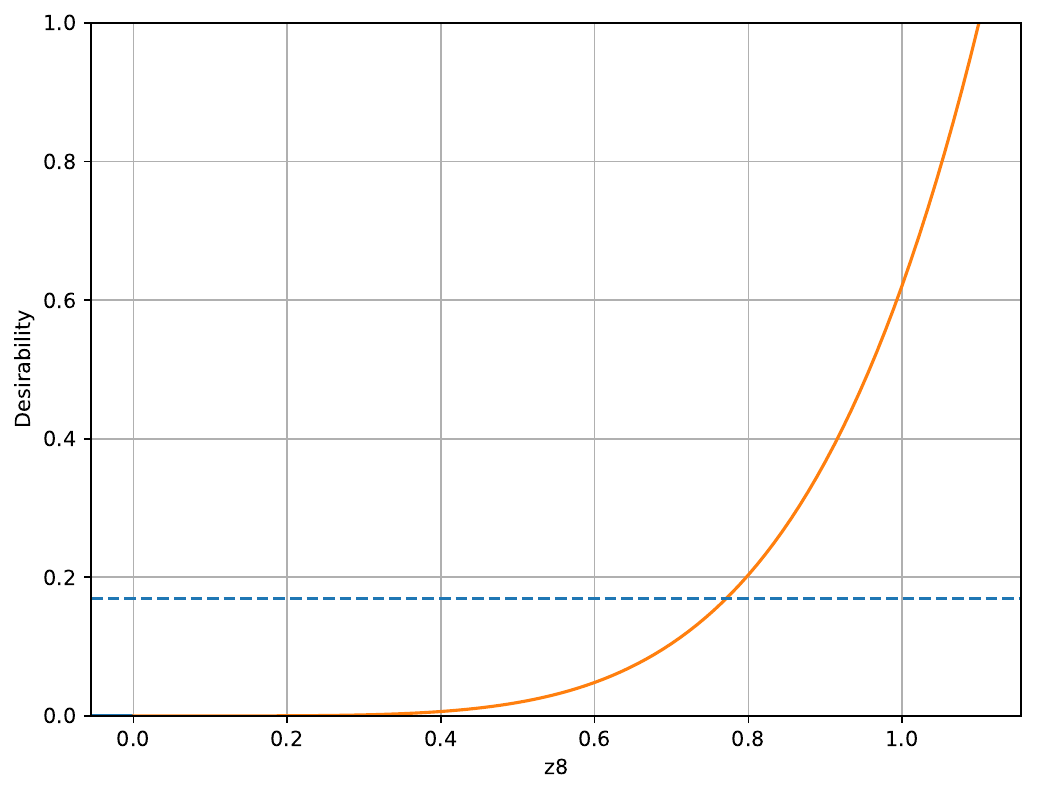}}
}
\caption{\label{fig-desirability-functions-z8-wo-mm}First desirability
function for target variable z8}
\end{figure}

\begin{figure}
\centering{
\pandocbounded{\includegraphics[width=0.8\linewidth,keepaspectratio]{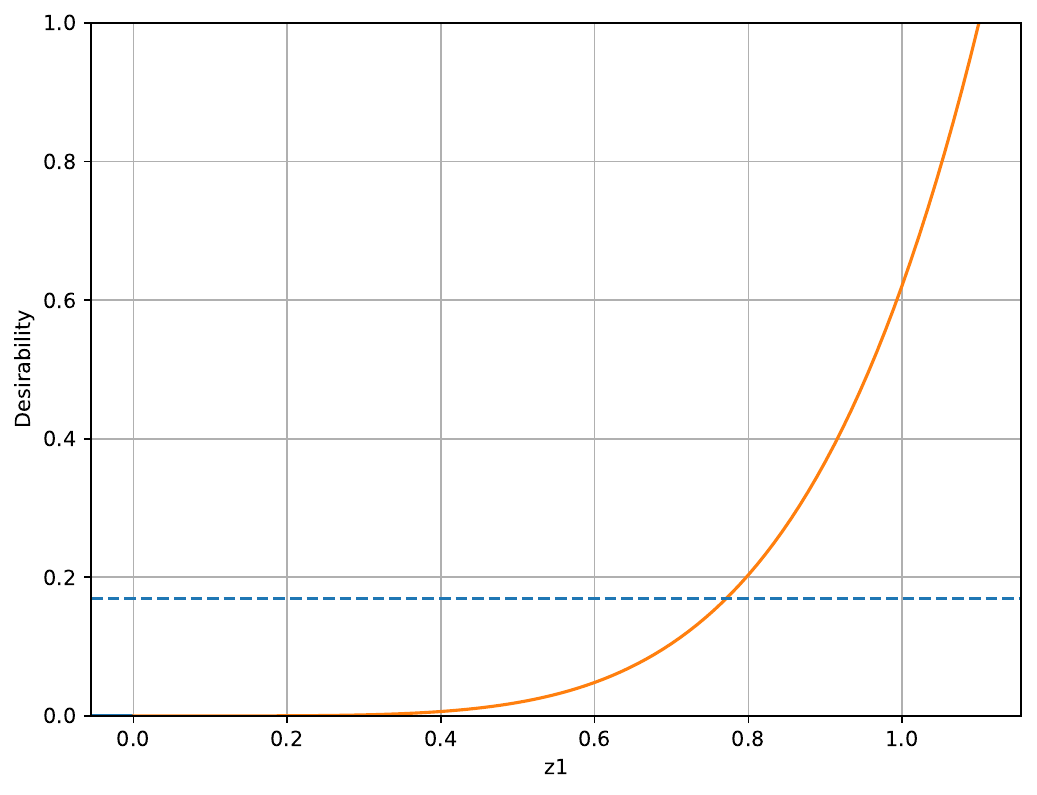}}
}
\caption{\label{fig-desirability-functions-z1-wo-mm}Second desirability
function for target variable z1}

\end{figure}

We have two target variables, \texttt{z8} and \texttt{z1}, each with its
own desirability function. Desirabilities for maximization can be
parameterized via the ``support'', i.e., via the interval in which the
desirability is greater than zero. This is given in
Equation~\ref{eq-Dmax} by the parameters $A$ and $B$. And they can
be parameterized via the ``steepness'', i.e., via the parameter $s$ in
Equation~\ref{eq-Dmax}. Here we have chosen $[5, 5]$, which pushes the
desirability functions to be more steep. The desirability functions are
visualized in Figure~\ref{fig-desirability-functions-z8-wo-mm} and
Figure~\ref{fig-desirability-functions-z1-wo-mm}.

The optimizer uses the desirability \texttt{z} as a combination of the
two objective functions for the search, i.e.,
$z = D_{\text{overall}}(f_1(x), f_2(x))$. Using the \texttt{DOverall}
function from the \texttt{spotdesirability} package, the two
desirabilities are combined into a single overall desirability. The
combined desirability is then used in the optimization process.

\subsection{Using the Combined Desirability in Multi-Objective
Optimization}\label{using-the-combined-desirability-in-multi-objective-optimization}

\begin{figure}
\centering{
\pandocbounded{\includegraphics[keepaspectratio]{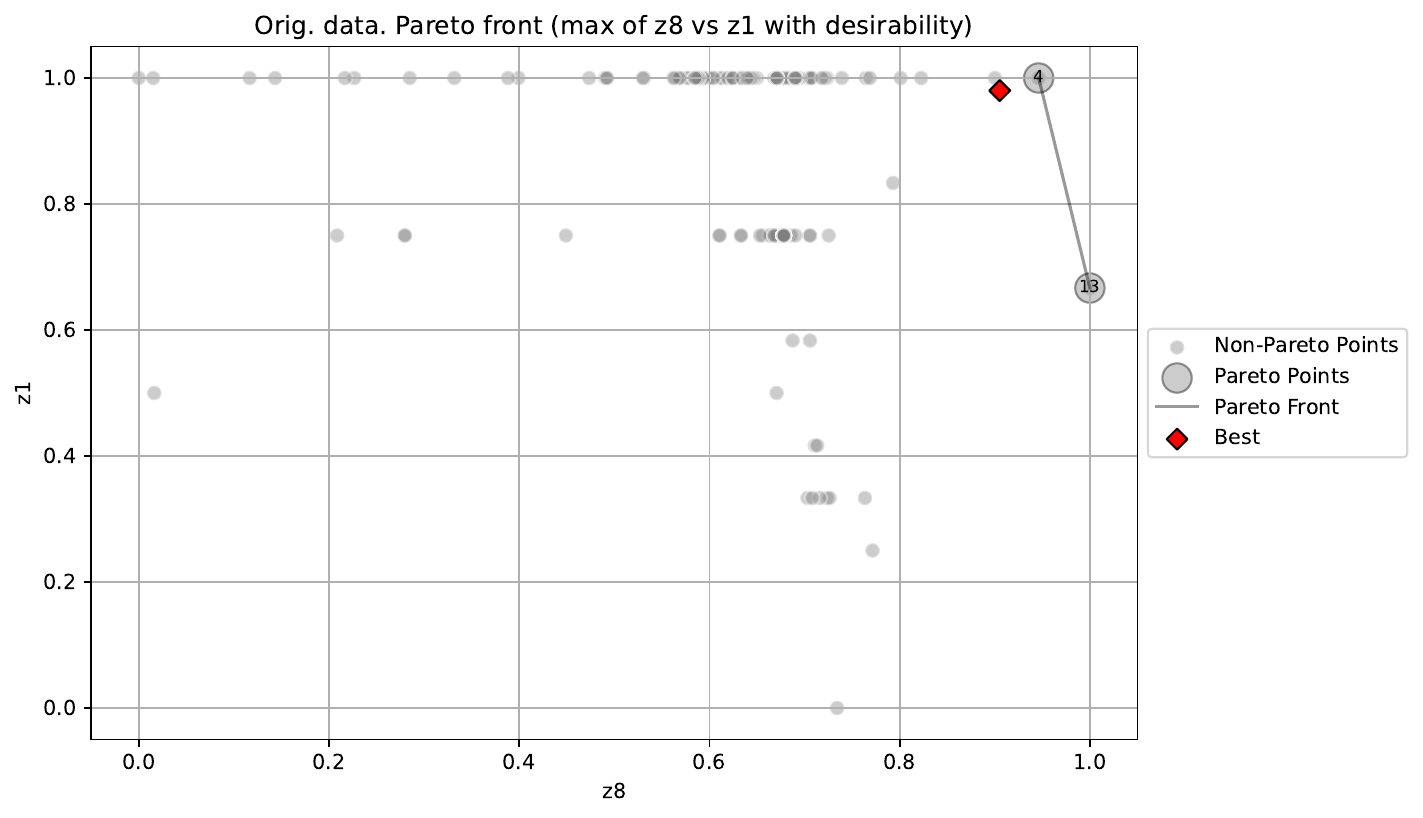}}
}
\caption{\label{fig-paretofront-orig-opt-z8-z1-zeo}Pareto front of the
original data and the optimized data for the multi-dimensional case (z8
vs.~z1)}
\end{figure}

\begin{figure}
\centering{
\pandocbounded{\includegraphics[keepaspectratio]{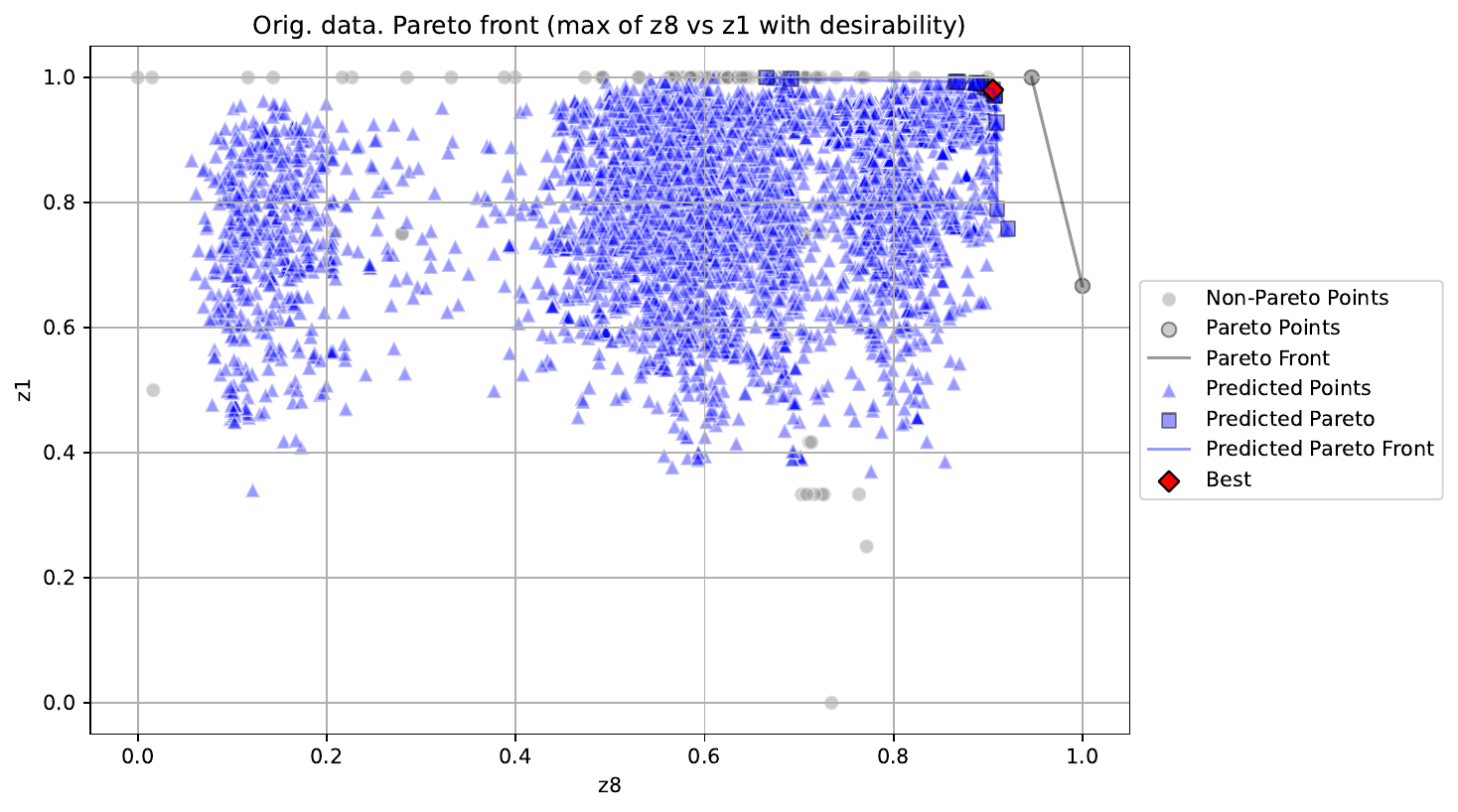}}
}
\caption{\label{fig-paretofront-orig-opt-z8-z1-zeo-callback}Pareto front
of the original data and the optimized data for the multi-dimensional
case (z8 vs.~z1)}

\end{figure}

In Section~\ref{sec-mm-criterion-multi}, we will use the intensified
Morris-Mitchell criterion $\Phi^{\ast}$ as an additional objective in
multi-objective optimization. Here, we first consider the case without
$\Phi^{\ast}$, which allows us a comparison with the case including
the Morris-Mitchell criterion. That is, we optimize the combined
desirability of the two target variables \texttt{z8} and \texttt{z1}.
The argument \texttt{mm\_objective=False} indicates that the
Morris-Mitchell criterion is not considered as an additional objective.
After running the multi-objecive optimization, we can calculate the best
target function values for \texttt{z8} and \texttt{z1} (here:
\texttt{y\_best\_z8\_z1}), i.e., the best predicted values from
multiobjective optimization of \texttt{z8} and \texttt{z1}, the best
input values for \texttt{z8} and \texttt{z1} (here:
\texttt{best\_x\_z8\_z1}), and the best desirability (here:
\texttt{best\_desirability\_z8\_z1}), i.e., the best desirability from
multiobjective optimization of \texttt{z8} and \texttt{z1}.
Figure~\ref{fig-paretofront-orig-opt-z8-z1-zeo} shows the Pareto front
of the original data and the optimized point \texttt{y\_best\_z8\_z1}
(shown in ``red'') for the multi-dimensional case (\texttt{z8} versus
\texttt{z1}). We can also add the \texttt{callback\_values}, i.e., the
values of the desirability function at each iteration of the
optimization, to the plot. This is shown in
Figure~\ref{fig-paretofront-orig-opt-z8-z1-zeo-callback}.

\section{\texorpdfstring{The Intensified Morris-Mitchell Criterion
$\Phi^{\ast}$ as an Additional
Objective}{The Intensified Morris-Mitchell Criterion \textbackslash Phi\^{}\{\textbackslash ast\} as an Additional Objective}}\label{sec-mm-criterion-multi}

By including the improvement of the search space coverage, the reduction
of the value of the intensified Morris-Mitchell criterion
$\Phi^{\ast}$ is added as the $p+1$st target variable to existing
multi-objective optimization problems with $p$ target variables, which
should be maximized. The intensified Morris-Mitchell criterion is used
here, because it shows similar results as the corrected Morris-Mitchell
criterion, as shown in Figure~\ref{fig-mm-corrected-vs-n-lhs}.
Furthermore, $\Phi^{\ast}$ is computationally less expensive than
$\hat{\Phi}$.

Let \begin{equation}\protect\phantomsection\label{eq-mm-improvement}{
f_3(x) = \Phi^{\ast}(X) - \Phi^{\ast}(X^+), \text{ where } X^+ = X \cup \{x\} \text{ and } X \text{ is the current search space.}
}\end{equation} This is in our case ($p=2$) the third objective
function to be maximized, and it represents the reduction in
$\Phi^{\ast}$ when adding the new point $x$ to the existing design
$X_u$. A positive value indicates that the new point improves the
space-filling quality of the design. Using a desirability function, the
three criteria are combined to obtain an overall desirability.

\subsection{\texorpdfstring{Determination of Desirability for
$\Phi^{\ast}$}{Determination of Desirability for \textbackslash Phi\^{}\{\textbackslash ast\}}}\label{determination-of-desirability-for-phiast}

We start with the calculation of the quality of the current design
$X_u$, which contains the entire dataset \texttt{X}, using
$\Phi^{\ast}$. The value of $\Phi^{\ast}(X_u)$ is used as a baseline
to evaluate the quality of the design.

\protect\phantomsection\label{mmphi-intensive-x-1krit}
\begin{verbatim}
mmphi_base (known design): 136.33458506472726
\end{verbatim}

If another point is added, the $\Phi^{\ast}$ value should be reduced
so that the difference between the Morris-Mitchell value of the current
design and the design with the added point, i.e., $f_3(x)$ as defined
in Equation~\ref{eq-mm-improvement}, should be positive and should be
maximized. In our case, we have chosen a relative reduction between
0.1\% and 2.5\%. The minimal and maximal values for the desirability
function of $\Phi^{\ast}$ are therefore set as follows:

\protect\phantomsection\label{mmphi-min-max-setup}
\begin{verbatim}
mmphi_min: 0.14, mmphi_max: 3.41
\end{verbatim}

Next, we combine the minimal and maximal $\Phi^{\ast}$ values with the
minimal and maximal values for the $z_8$ and $z_1$ objectives and
generate, as above, the overall-desirability function\footnote{Note, to
  implement this, we use the scale parameters
  \texttt{z8\_min\_multiplier=1.0} and \texttt{z8\_max\_multiplier=1.1}
  as well as \texttt{z1\_min\_multiplier} and
  \texttt{z1\_max\_multiplier=1.1} as before.}.

\protect\phantomsection\label{append-mmphi-min-max-234}
\begin{verbatim}
target z8: min: 0.0, max: 1.1, scale: 5
target z1: min: 0.0, max: 1.1, scale: 5
target mm: min: 0.13633458506472726, max: 3.408364626618182, scale: 5
\end{verbatim}

\begin{figure}
\centering{
\pandocbounded{\includegraphics[width=0.8\linewidth,keepaspectratio]{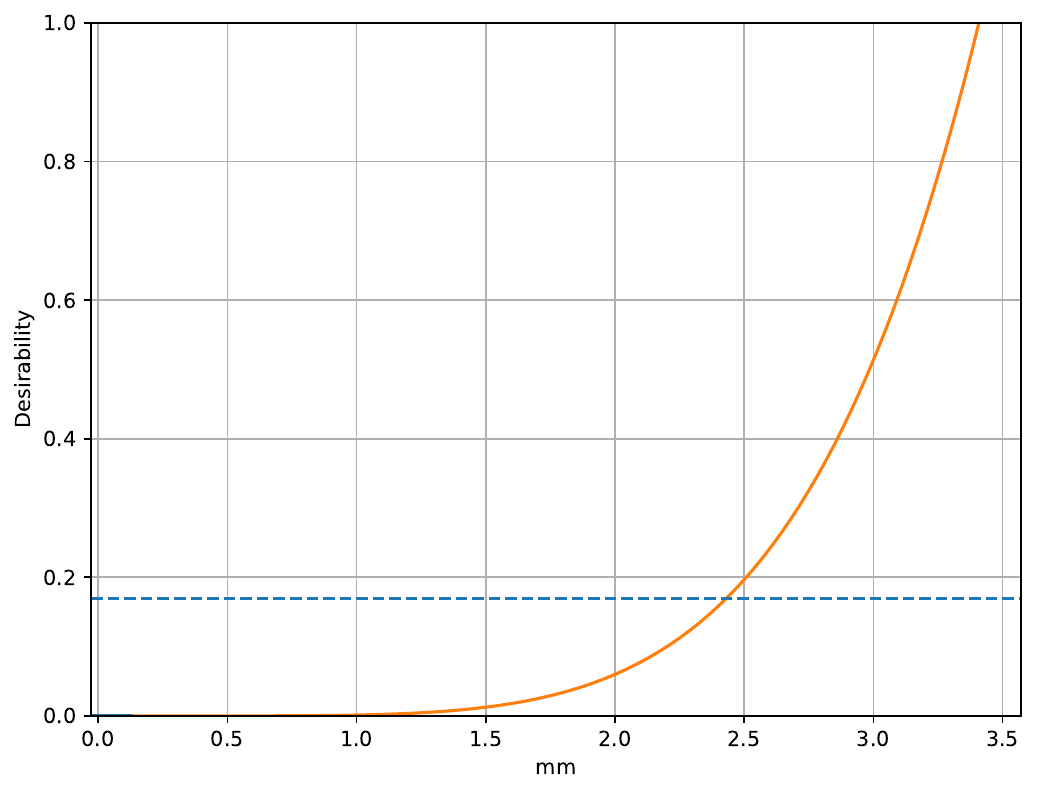}}
}
\caption{\label{fig-desirability-functions-mm-1}Third desirability
function (Morris-Mitchell criterion)}

\end{figure}

The desirability functions for the target variables \texttt{z8} and
\texttt{z1} were already shown in
Figure~\ref{fig-desirability-functions-z8-wo-mm} and
Figure~\ref{fig-desirability-functions-z1-wo-mm}, respectively.
Figure~\ref{fig-desirability-functions-mm-1} shows the third
desirability function for the Morris-Mitchell criterion.

\subsection{Optimization with Three
Objectives}\label{optimization-with-three-objectives}

Now we have three desirability functions that are combined into an
overall desirability function \texttt{overallD}. We can consider the
optimization including $\Phi^{\ast}$ as an additional objective.
Running optimization for \texttt{z8}, \texttt{z1} and $\Phi^{\ast}$
with \texttt{mo\_mm\_desirability\_function} and
\texttt{mm\_objective=True}. By setting \texttt{verbose=True} in the
function \texttt{mo\_mm\_desirability\_optimizer}, the $\Phi^{\ast}$
values during optimization can be printed.

\protect\phantomsection\label{verify-mo-function-z8-z1-mm}
\begin{verbatim}
Best point (z8 + z1 + MM): [0.97488619 0.53048699 0.58165396 0.91648553 0.86938054 0.98930395
 0.14399584 0.55196457 0.18361494 0.05696308 0.6921046  0.09899566
 0.08910941 0.79756473 0.85654689 0.38576957 0.66376137 0.45312888
 0.82637003 0.5310207  0.74631296 0.9347274  0.65727871 0.88256815
 0.73470113 0.92485178 0.97761449]
Best Desirability (z8 + z1 + MM): 0.2088
y_best (z8, z1, MM): [[0.92232323 0.95416667 1.89376208]]
\end{verbatim}

\begin{figure}
\centering{
\pandocbounded{\includegraphics[keepaspectratio]{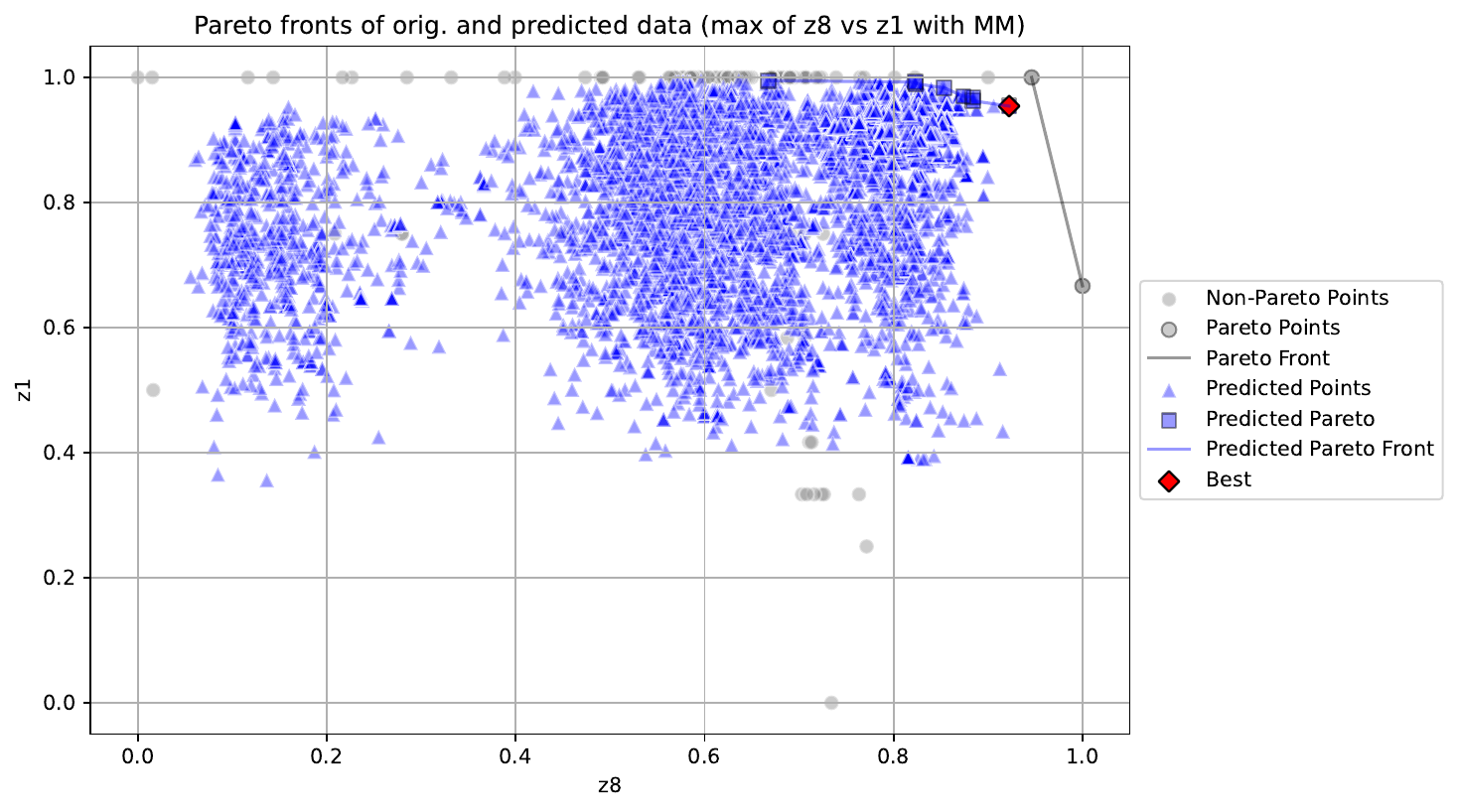}}
}
\caption{\label{fig-paretofront-orig-opt-z8-z1-mm}Pareto front of the
original data and the optimized data for the multi-dimensional case (z8
vs.~z1) with $\Phi^{\ast}$}
\end{figure}

\begin{figure}
\centering{
\pandocbounded{\includegraphics[keepaspectratio]{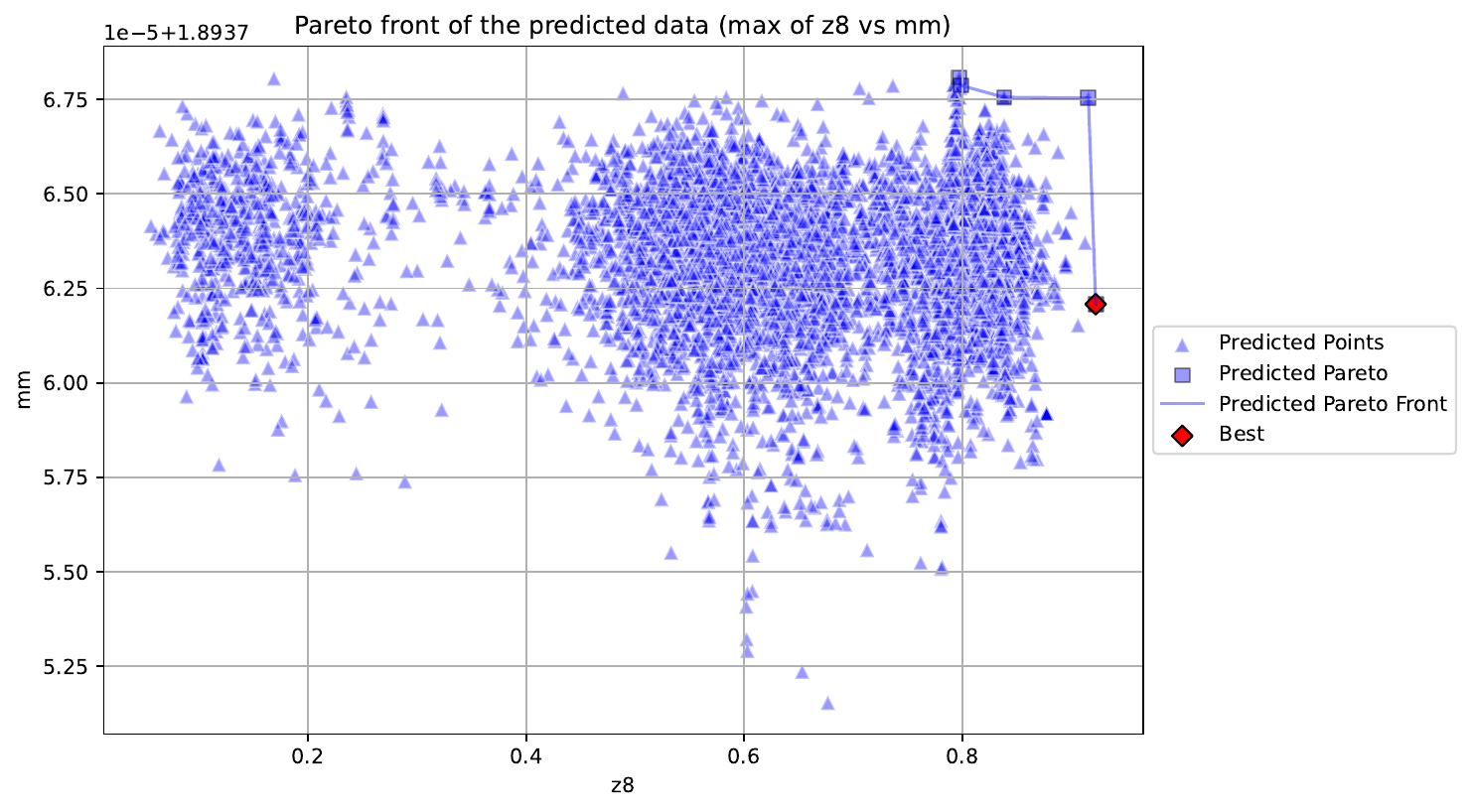}}
}
\caption{\label{fig-paretofront-orig-opt-z8-mm}Pareto front of the
optimized data for the multi-dimensional case (z8 vs.~$\Phi^{\ast}$).
For the original data, no $\Phi^{\ast}$ values are available.}

\end{figure}

\begin{figure}

\centering{

\pandocbounded{\includegraphics[keepaspectratio]{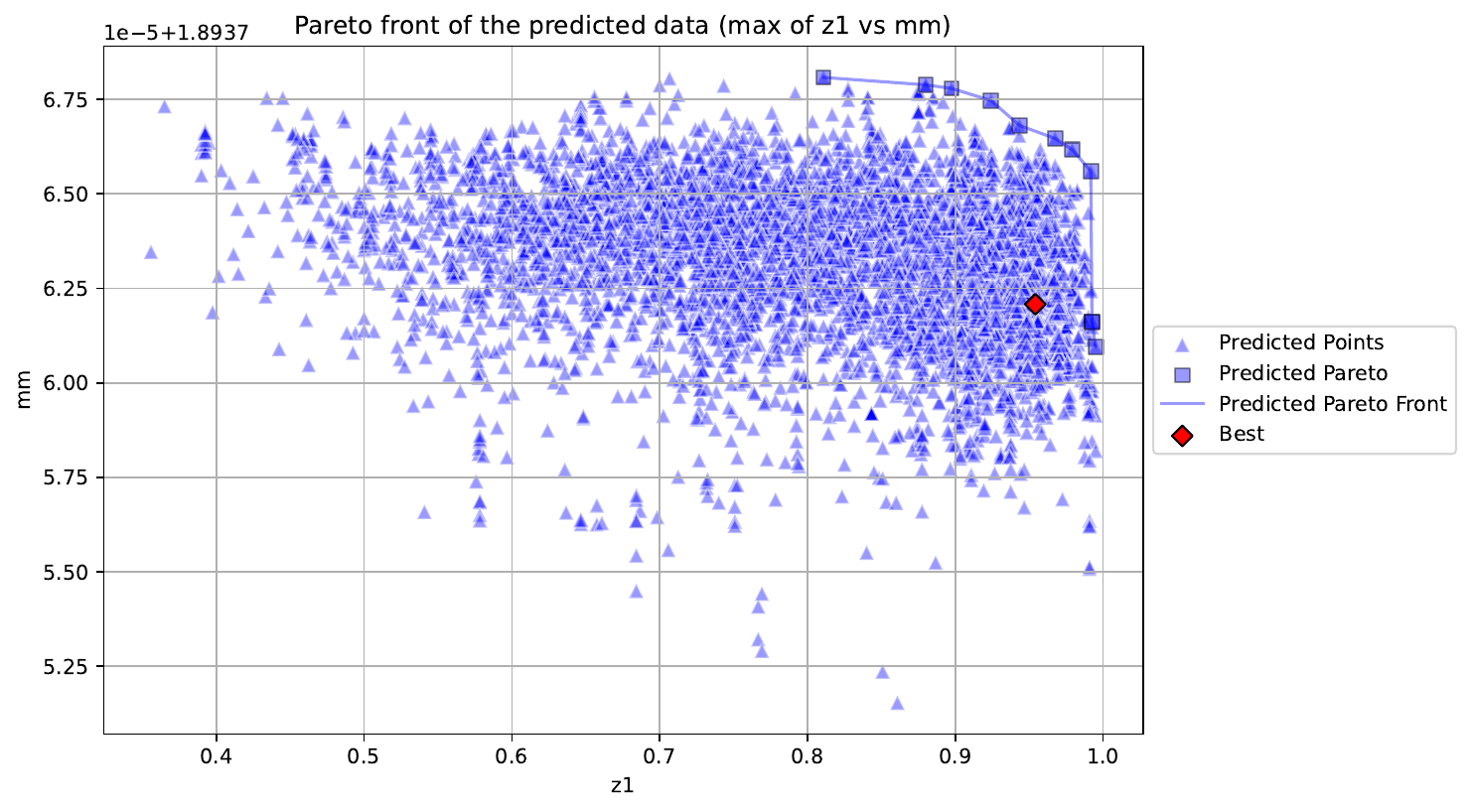}}

}

\caption{\label{fig-paretofront-orig-opt-z1-mm}Pareto front of the
optimized data for the multi-dimensional case (z1 vs.~$\Phi^{\ast}$).
For the original data, no mm values are available.}

\end{figure}

Similar to the results in Section~\ref{sec-opt-wo-mm}, we get the best
point for \texttt{z8} and \texttt{z1}, its desirability and the
corresponding target values for \texttt{z8} and \texttt{z1}. There are
three combinations of $\Phi^{\ast}$ with \texttt{z8} and \texttt{z1},
i.e., we can plot the Pareto front for \texttt{z8} vs.~\texttt{z1},
\texttt{z8} vs.~\texttt{mm}, and \texttt{z1} vs.~\texttt{mm}:
Figure~\ref{fig-paretofront-orig-opt-z8-z1-mm} shows the Pareto front of
the original data and the optimized data for the multi-dimensional case
(\texttt{z8} vs.~\texttt{z1}) with $\Phi^{\ast}$. The best point is
colored in ``red''. In the following two cases, no original data is
available, since $\Phi^{\ast}$ is only calculated for the optimized
points. Figure~\ref{fig-paretofront-orig-opt-z8-mm} shows the Pareto
front of the optimized data for the multi-dimensional case (\texttt{z8}
vs.~$\Phi^{\ast}$). Figure~\ref{fig-paretofront-orig-opt-z1-mm} shows
the Pareto front of the optimized data for the multi-dimensional case
(\texttt{z1} vs.~$\Phi^{\ast}$).

\section{Results and Discussion}\label{sec-results-discussion}

The two case studies show how to determine the next design point for a
given, unplanned design space $X_u$. A surrogate model was trained on
the existing data from the original design space $X_u$ and used to
predict promising infill points. This is a standard approach in
surrogate-model based optimization (Forrester et al. 2008). Similar to
the expected improvement criterion, which balances the trade-off between
exploration and exploitation, the Morris-Mitchell criterion can be used
to determine the next design point. The Morris-Mitchell criterion
balances the trade-off between good design points and the
space-fillingness of the updated design space $X_u \cup x$.

To investigate the impact of adding the improved Morris-Mitchell
criterion $\Phi^{\ast}$ to the optimization, we compare the Pareto
front of the original data and the optimized data for the
multi-dimensional case (\texttt{z8} vs.~\texttt{z1}) with
$\Phi^{\ast}$ to the Pareto front of the original data and the
optimized data for the multi-dimensional case (\texttt{z8}
vs.~\texttt{z1}) without $\Phi^{\ast}$. This comparison reveals no
significant difference in the optimization with and without
$\Phi^{\ast}$, which gives rise to two considerations.

\begin{itemize}
\tightlist
\item
  On the one hand, it should be noted that the underlying dataset was
  already extensively optimized manually and therefore already
  represents the Pareto front and optimum quite well.
\item
  On the other hand, it would be interesting to repeat this study with a
  significantly reduced dataset that does not contain the already
  existing optima. This would help answer whether the current Pareto
  front and optimum would still be found under such conditions.
\end{itemize}

A benefit of adding $\Phi^{\ast}$ to the classical multi-objective
optimization is that it considers the space-fillingness of the updated
design space $X_u \cup x$.\\
We claim that the results of this optimization can give the practitioner
valuable information where to place the next design point. As an
additional diagnostic, we generate infill-point histograms
(ip-histograms) and infill-point boxplots (ip-boxplots) to visualize the
distribution of the existing data points and the newly suggested best
point in the context of the existing data.

\subsection{Infill-point Diagnostics}\label{sec-ip-diagnostics}

Infill-point diagnostic plots are comprehensive tools to visualize the
location of the newly suggested best point in the context of the
existing data (Bartz-Beielstein 2025a). In the following figures, the
best points from the optimization with (\texttt{red}) and without
(\texttt{blue}) $\Phi^{\ast}$ are visualized.

\begin{figure}

\centering{

\pandocbounded{\includegraphics[keepaspectratio]{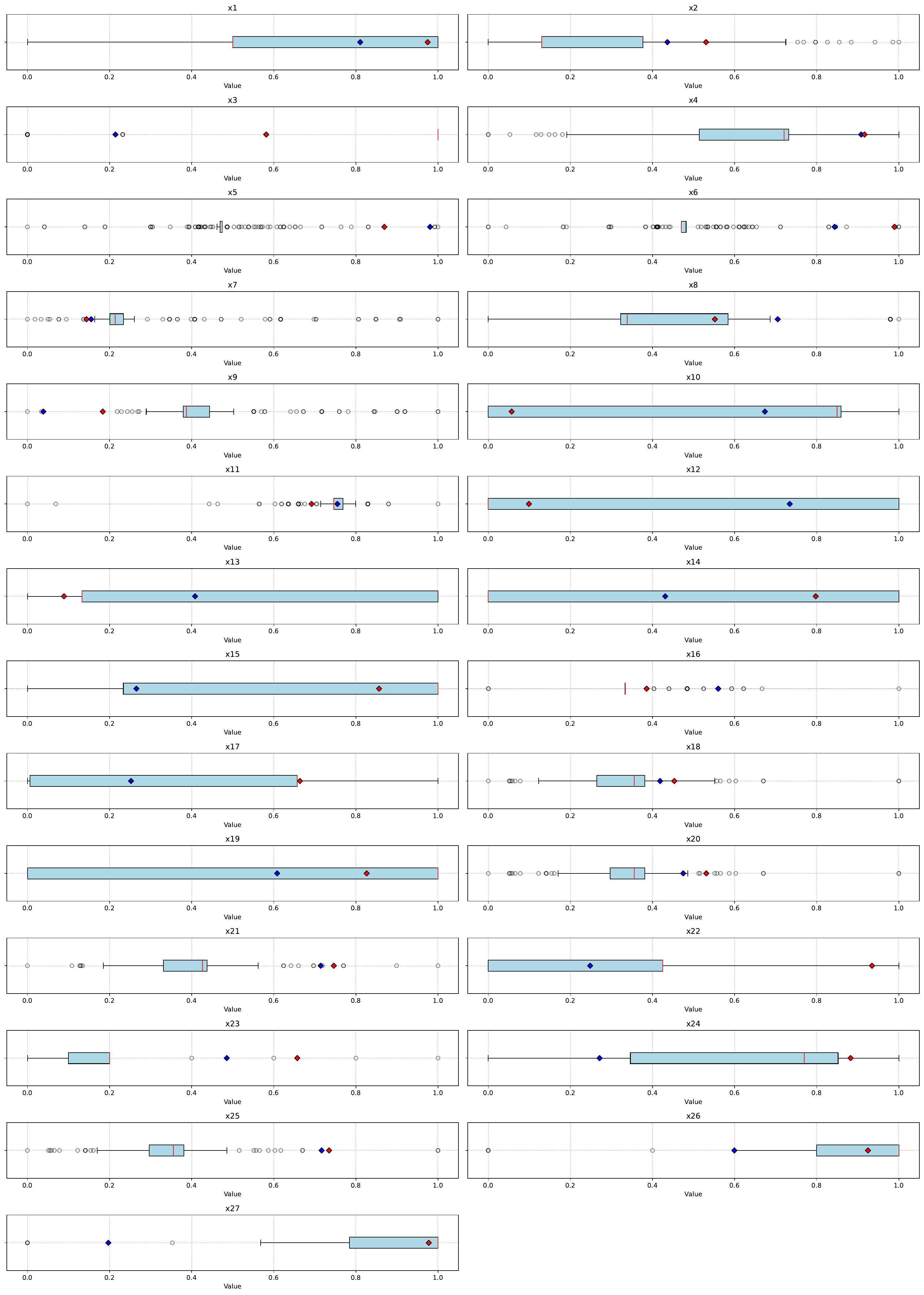}}

}

\caption{\label{fig-ip-box-1}Boxplots with the new setting with
$\Phi^{\ast}$ (red) and without $\Phi^{\ast}$ (blue)}

\end{figure}

\begin{figure}

\centering{

\pandocbounded{\includegraphics[keepaspectratio]{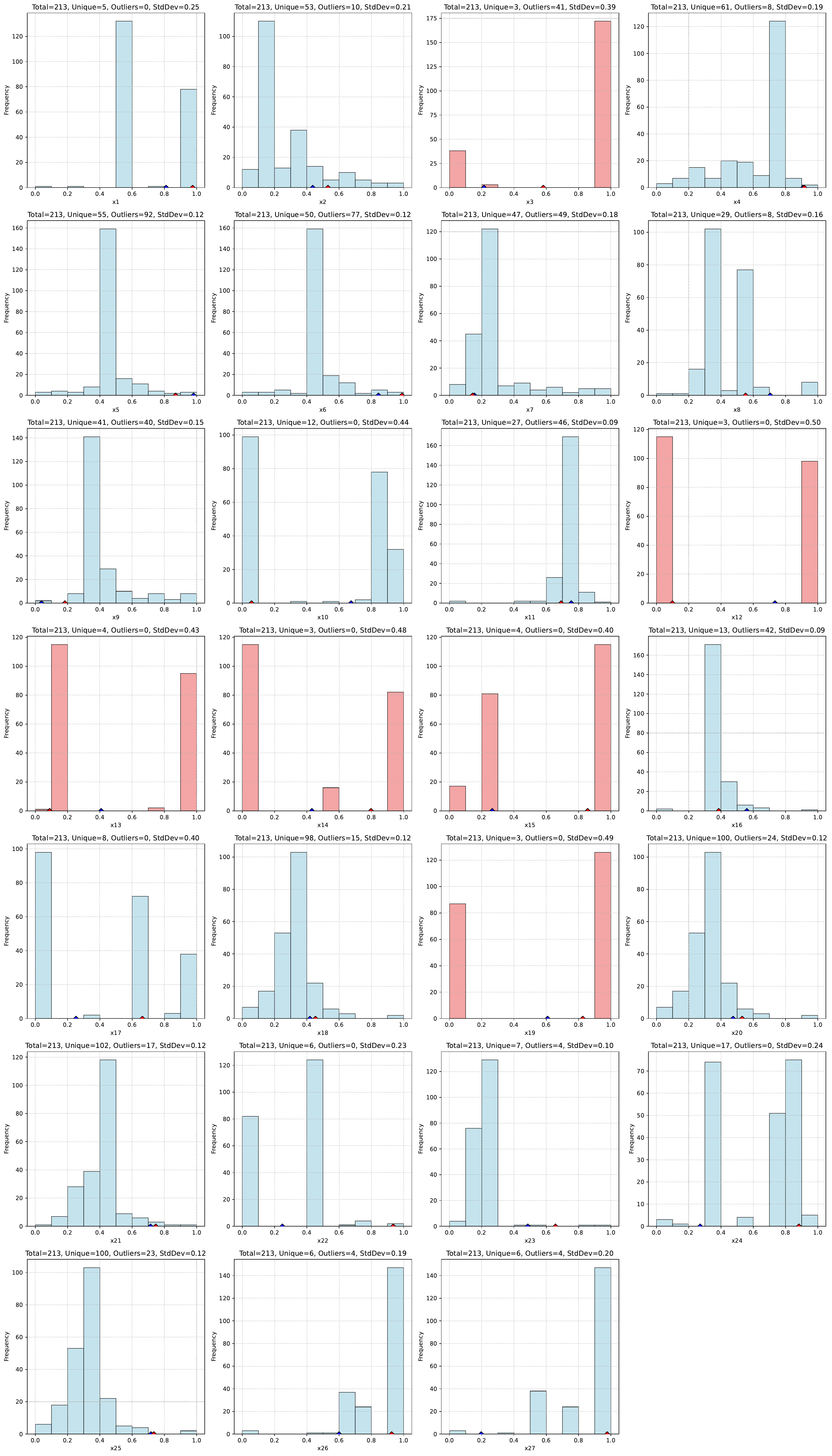}}

}

\caption{\label{fig-ip-hist-1}Histograms with the new setting with
$\Phi^{\ast}$ (red) and without $\Phi^{\ast}$ (blue). Histogram
plotted in red indicate a small number of distinct values.}

\end{figure}

Considering the first input variable, $x_1$, in
Figure~\ref{fig-ip-box-1}, we can see that an increase of its value is
recommended. Including the $\Phi^{\ast}$, the best value is
$x_1 \approx 1.0$, whereas the optimization without $\Phi^{\ast}$
leads to $x_1 \approx 0.8$. For $x_2$, both optimizations recommend
an increase ($x_2 \approx 0.5$). Figure~\ref{fig-ip-hist-1} uses the
same data as Figure~\ref{fig-ip-box-1} but shows histograms instead of
boxplots.

\begin{figure}

\centering{

\pandocbounded{\includegraphics[keepaspectratio]{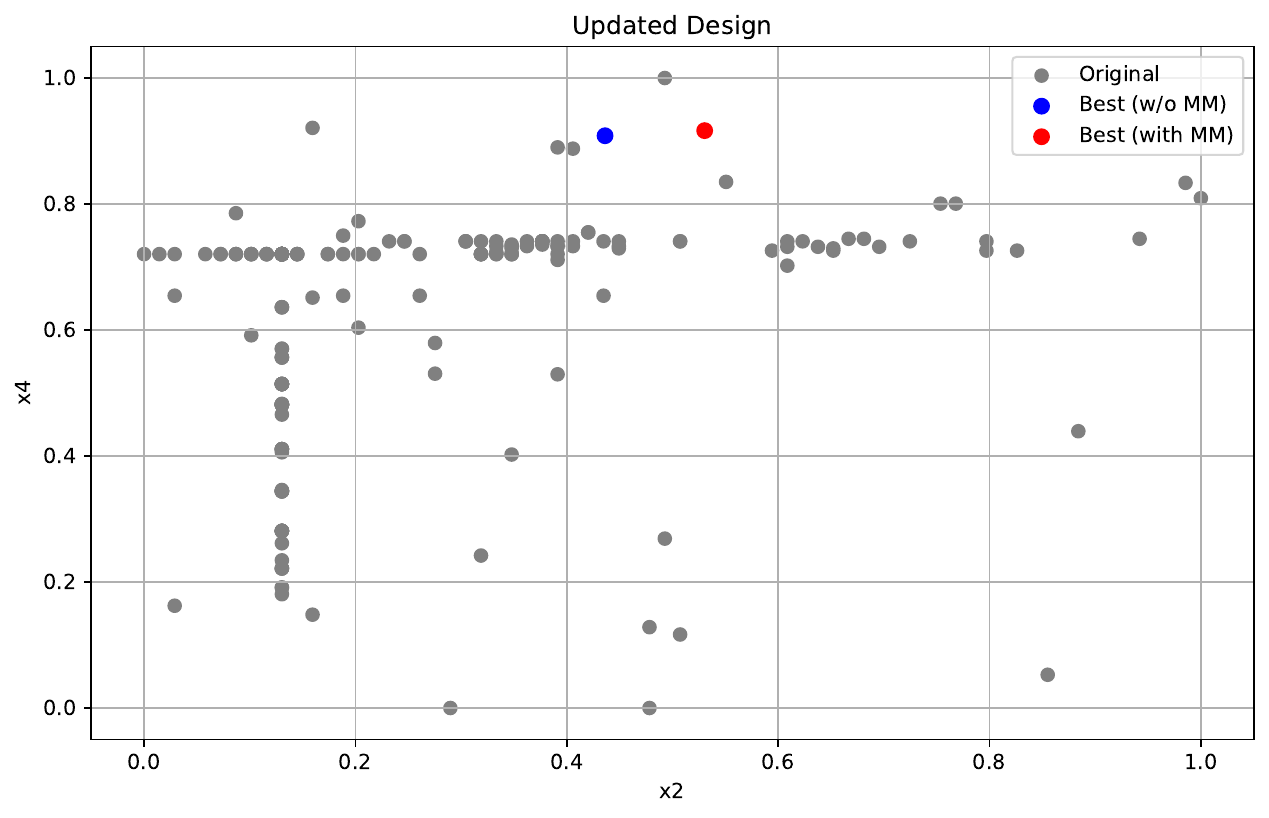}}

}

\caption{\label{fig-man-design-updated}Compressor Data. Updated Design.
Original data is marked in grey. The newly suggested best point using
the $\Phi^\ast$ criterion is marked in red. The best point without
$\Phi^\ast$ is marked in blue.}

\end{figure}

Finally, we reconsider the initial motivation of this study: the
clustered design space $X_u$ as shown in Figure~\ref{fig-man-design}.
Figure~\ref{fig-man-design-updated} shows the updated design space
$X_u \cup x$, i.e., the location of the newly suggested best point in
the context of the existing data.

\subsection{Improvement of the Intensified Morris-Mitchell Criterion by
Adding
Points}\label{improvement-of-the-intensified-morris-mitchell-criterion-by-adding-points}

Comparing the results from the two studies (optimization with
Morris-Mitchell criterion and without Morris-Mitchell criterion), one
might conclude that the Morris-Mitchell criterion has only a minor
effect on the outcome. For a further analysis, we consider in this
section the intensified Morris-Mitchell criterion for different numbers
of added points. Figure~\ref{fig-mmphi-intens-plot-x-call} shows the
values of the intensified Morris-Mitchell criterion, $\Phi^{\ast}$
versus the number of added points. We can see from this figure, that by
adding ten points to the existing design, the Morris-Mitchell criterion
reduced from approximatively 136 to 130, indicating an improvement in
the space-filling property of the design.

\protect\phantomsection\label{mmphi-intensive-x}
\begin{verbatim}
mmphi_base (known design): 136.33458506472726
mmphi_base extended with 10 random points: 130.20722322230748
\end{verbatim}

This reduction in the $\Phi^{\ast}$ value motivated the usage of the
intensified Morris-Mitchell criterion as an additional objective in the
optimization. If this reduction is large, a better design is found. But
how large can this reduction be? Is it a good indicator for an
multi-objective optimization algorithm?

To analyse these questions, we consider updating $\Phi^{\ast}$ by
adding only one point. This is motivated by the fact that in practice,
we usually add one point at a time in the optimization process, which
will be referred to as ``single-point injection''. Note, not the same
point is added several times, but different points, each randomly
generated. This allows us to determine the variance of the improvement
of $\Phi^{\ast}$.

Adding one random point to the existing compressor design reduces the
$\Phi^{\ast}$ by approximatively 0.64, no matter, which point is
injected. This indicates that adding any random point improves the
$\Phi^{\ast}$ value by a similar value. There is nearly no variance in
adding different random points, because all random points approach the
maximum possible improvement in the intensified Morris-Mitchell
criterion.

The reason is that the existing points are not very well distributed in
the design space, so that adding \emph{any} random point improves the
space-filling property of the design significantly.

\begin{figure}
\centering{
\pandocbounded{\includegraphics[width=0.8\linewidth,keepaspectratio]{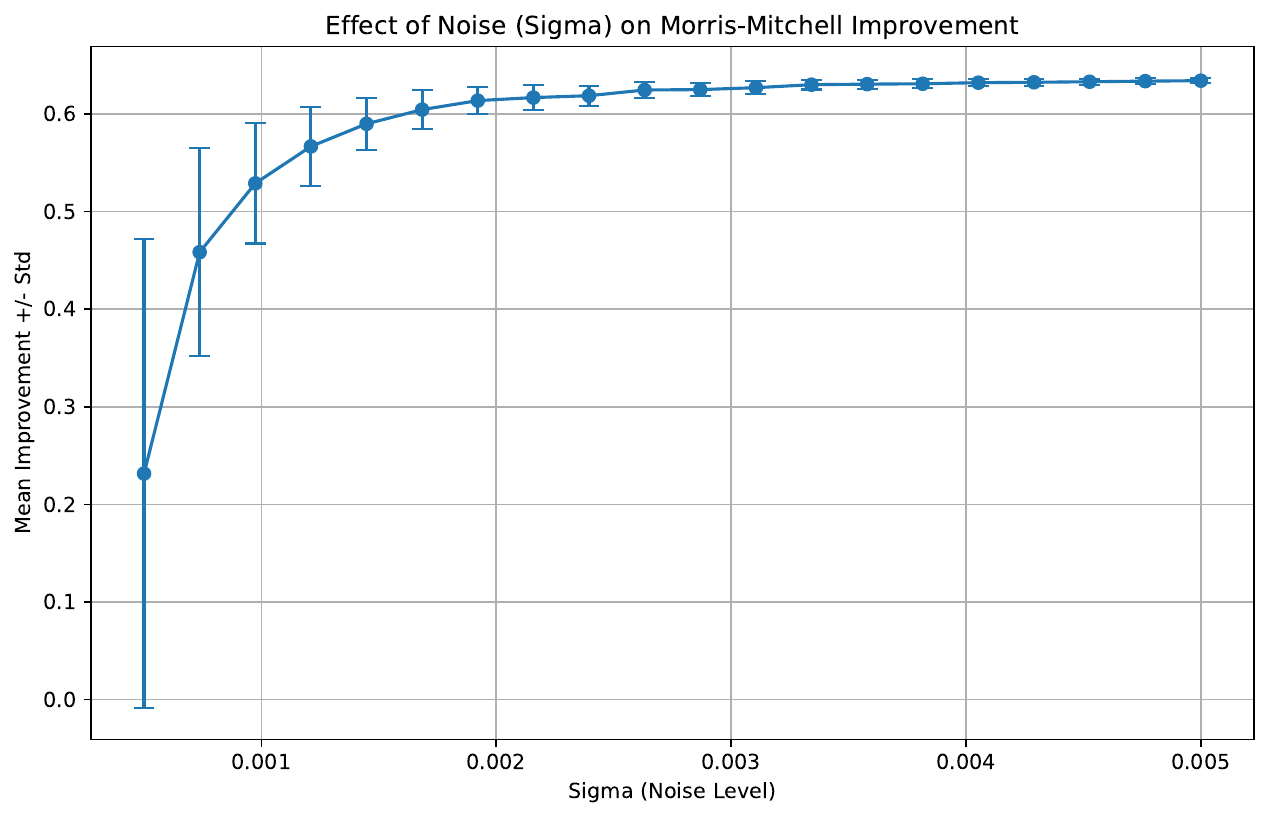}}
}
\caption{\label{fig-sigma-noise-1}Compressor Dataset. Effect of Noise
(Sigma) on the intensified Morris-Mitchell criterion improvement. Higher
noise (sigma) values generate infill points that are further away from
existing points.}

\end{figure}

To analyze this effect further, we inject points that are close to
existing points. Now we consider updating the Morris-Mitchell criterion
by adding ten times one point that is close to existing points instead
of adding random points that are not close to the existing design. This
allows us to see the effect of adding each point individually on the
Morris-Mitchell criterion. We add small noise to the existing points to
create new points that are close to the existing points, but its
distance to the existing points is modified by the noise added. Then we
compute the mean improvement of $\Phi^{\ast}$ based on these ten
``single-point injections'' and the standard deviation of the
improvement for each noise level (sigma). Results are shown in
Figure~\ref{fig-sigma-noise-1}. This figure reveals some interesting
insights. Points close to the existing design improve $\Phi^{\ast}$
only by a small amount, but the variance is higher than for random
points. By moving further away from the existing design by increasing
the noise level, the improvement increases and reaches its maximum at
0.64. This value is identical to the improvement of randomly added
points as discussed above.

These empirical findings can be directly connected to the theoretical
results derived in this paper. From the proof of
Theorem~\ref{thm-sensitivity-mmphi}, the improvement in $\Phi_q^{*q}$
when adding $x_{n+1}$ is exactly \[
\Phi_q^{*q}(X_p) - \Phi_q^{*q}(X_p \cup \{x_{n+1}\}) = \frac{2}{n+1}\left[\Phi_q^{*q}(X_p) - \frac{\Delta(x_{n+1}, X_p)}{n}\right],
\] where $\Delta(x_{n+1}, X_p) = \sum_{k=1}^{m'} J_k'(d_k')^{-q}$ is
the interaction energy of the new point with the existing design. For
$n = 213$, the prefactor $2/(n+1) \approx 0.0093$ is essentially
constant, so the improvement depends entirely on
$\Phi_q^{*q} - \Delta/n$. When the existing design is bad
($\Phi_q^{*q}$ large) and the new point is placed far from clusters
($\Delta/n$ moderate), this difference is large and stable ---
explaining the constant improvement of $\approx 0.64$. When a point is
injected close to an existing point (small sigma), the $d^{-q}$ terms
in $\Delta$ explode for the nearest neighbors, driving $\Delta/n$ up
toward $\Phi_q^{*q}$ and reducing the improvement. This is exactly the
pattern in Figure~\ref{fig-sigma-noise-1}.

\subsection{Summary and Outlook}\label{sec-summary-outlook}

This study presents a comprehensive methodology for optimizing
experimental designs through a surrogate-model-based approach that
integrates desirability functions and the Morris-Mitchell criterion. The
application of multi-objective optimization utilizing desirability
functions yields transparent and intuitive results that effectively
support practical decision-making. To quantify and improve the quality
of prospective infill points, we rigorously introduce the intensified
Morris-Mitchell criterion, supported by detailed theoretical derivations
that elucidate the intrinsic relationship between the intensified and
the corrected Morris-Mitchell criteria.

Furthermore, the integration of potential theory provides a robust
mathematical framework for analyzing space-filling geometries. From this
perspective, minimizing the Morris-Mitchell criterion within statistical
experimental design is fundamentally equivalent to identifying a point
configuration that minimizes its internal Riesz energy, representing the
total repulsive force among all points. This energy minimization
naturally forces points to separate maximally, thereby achieving a
uniform spatial distribution. We leverage this structural connection to
critically analyze the monotonicity properties of the Morris-Mitchell
criterion. Finally, to facilitate the practical realization of these
theoretical advancements, we propose novel infill-point diagnostics.
These visualization tools, employing specialized boxplots and
histograms, offer practitioners an effective, geometry-aware guide for
optimally placing subsequent design points within previously unplanned
experimental spaces.

These modern theoretical and diagnostic techniques are practically
implemented and accessible via the Python packages
\texttt{spotdesirability} and \texttt{spotoptim}.
\texttt{spotdesirability}\footnote{Documented at
  \url{https://sequential-parameter-optimization.github.io/spotdesirability/}}
implements the desirability framework outlined by Derringer \& Suich
(1980), providing a versatile set of professional tools for multifaceted
multi-objective optimization. It serves as a Python port of the esteemed
R package \texttt{desirability}\footnote{Available on CRAN:
  \url{https://CRAN.R-project.org/package=desirability}, DOI:
  \href{https://doi.org/10.32614/CRAN.package.desirability}{10.32614/CRAN.package.desirability}},
which was developed and modeled by Kuhn around the same S3 class
architecture. Likewise, \texttt{spotoptim}\footnote{Documented at
  \url{https://sequential-parameter-optimization.github.io/spotoptim/docs/index.html}}
serves as the primary optimizer for this workflow, integrating these
surrogate models seamlessly.

We will consider the following directions to further develop the
approach presented in this paper.
The Morris-Mitchell criterion, which we have treated as a geometric
measure and as a discrete Riesz energy functional, admits a statistical
interpretation through its connection to Kriging. In the Kriging
framework, the correlation between responses at two design points is
modeled as a function of their spatial distance via a stationary
Gaussian process (Antognini and Zagoraiou 2010).
A fundamental limitation of maximin-distance criteria, such as the
Morris-Mitchell criterion used in this paper, is the inherent lack of
strict monotonicity. As was discussed in Section~\ref{sec-MM-criterion},
in a bounded space, the minimum distance between points structurally
decreases as the number of points grows, preventing maximin-distance
criteria from monotonically decreasing when new points are added.
Instead, the dual perspective of covering quality offers natural
monotonicity properties. In particular, the covering radius (minimax
distance, the classical counterpart to the maximin criterion or
Morris-Mitchell criterion used in this paper) might be a promising
alternative (Pronzato 2017). Coverage criteria strictly decrease upon
any point addition, which would allow for easier integration and usage
of monotonicity properties in the optimization process. Furthermore,
this paradigm shift avoids needing normalizations by the number of
points or dimensionality.
Last but not least, we will consider some numerical aspects of the
presented approach. To avoid problems with local optima, it is useful to
perform restarts. Although we did not use restarts in our experiment, we
recommend its generation in the context of extending unplanned
industrial designs $X_u$. The set of restart points can be generated
using LHS in combiation with the Morris-Mitchell criterion.
Hyperparameter tuning should be applied to the model to improve the
quality of the predictions. \texttt{spotoptim} provides inherently a
good starting point for hyperparameter tuning.

\section{Appendix}\label{appendix}

\subsection{Cross-Validated Comparison of the
Models}\label{sec-cv-model-comparison}

\begin{figure}

\centering{

\pandocbounded{\includegraphics[keepaspectratio]{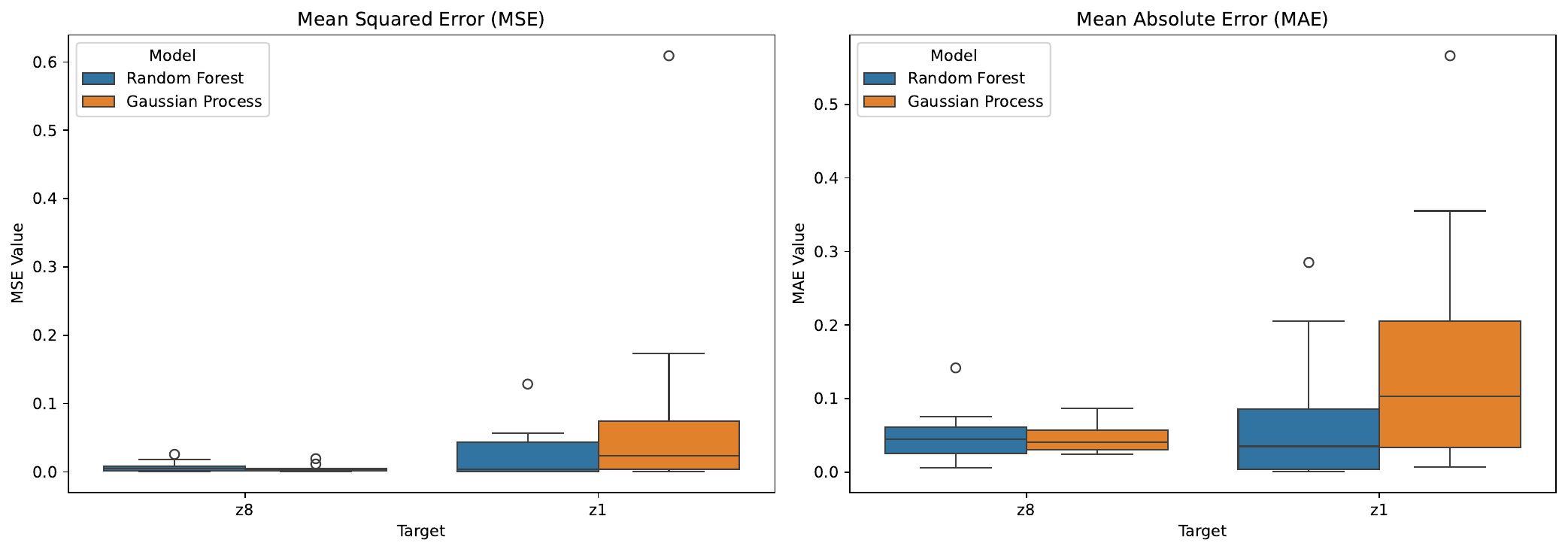}}

}

\caption{\label{fig-cv-plot-comparison}Cross-validated comparison of
Random Forest and Gaussian Process models (MSE and MAE). Smaller values
are better.}

\end{figure}

Table~\ref{tbl-cv-comparison} shows the cross-validated comparison of
the Random Forest and Gaussian Process models for the two target
variables \texttt{z8} and \texttt{z1}. The table includes the mean,
standard deviation, minimum and maximum values of the Mean Squared Error
(MSE) and Mean Absolute Error (MAE) across the 10 folds of
cross-validation. These results are visualized in
Figure~\ref{fig-cv-plot-comparison}.

\begin{longtable}[]{@{}
  >{\raggedright\arraybackslash}p{(\linewidth - 12\tabcolsep) * \real{0.1429}}
  >{\raggedright\arraybackslash}p{(\linewidth - 12\tabcolsep) * \real{0.2571}}
  >{\raggedright\arraybackslash}p{(\linewidth - 12\tabcolsep) * \real{0.1429}}
  >{\raggedleft\arraybackslash}p{(\linewidth - 12\tabcolsep) * \real{0.1143}}
  >{\raggedleft\arraybackslash}p{(\linewidth - 12\tabcolsep) * \real{0.1143}}
  >{\raggedleft\arraybackslash}p{(\linewidth - 12\tabcolsep) * \real{0.1143}}
  >{\raggedleft\arraybackslash}p{(\linewidth - 12\tabcolsep) * \real{0.1143}}@{}}
\caption{Cross-validated comparison of model performance (MSE and MAE).
Smaller values are better.}\label{tbl-cv-comparison}\tabularnewline
\toprule\noalign{}
\begin{minipage}[b]{\linewidth}\raggedright
Target
\end{minipage} & \begin{minipage}[b]{\linewidth}\raggedright
Model
\end{minipage} & \begin{minipage}[b]{\linewidth}\raggedright
Metric
\end{minipage} & \begin{minipage}[b]{\linewidth}\raggedleft
Mean
\end{minipage} & \begin{minipage}[b]{\linewidth}\raggedleft
Std
\end{minipage} & \begin{minipage}[b]{\linewidth}\raggedleft
Min
\end{minipage} & \begin{minipage}[b]{\linewidth}\raggedleft
Max
\end{minipage} \\
\midrule\noalign{}
\endfirsthead
\toprule\noalign{}
\begin{minipage}[b]{\linewidth}\raggedright
Target
\end{minipage} & \begin{minipage}[b]{\linewidth}\raggedright
Model
\end{minipage} & \begin{minipage}[b]{\linewidth}\raggedright
Metric
\end{minipage} & \begin{minipage}[b]{\linewidth}\raggedleft
Mean
\end{minipage} & \begin{minipage}[b]{\linewidth}\raggedleft
Std
\end{minipage} & \begin{minipage}[b]{\linewidth}\raggedleft
Min
\end{minipage} & \begin{minipage}[b]{\linewidth}\raggedleft
Max
\end{minipage} \\
\midrule\noalign{}
\endhead
\bottomrule\noalign{}
\endlastfoot
z8 & Random Forest & MSE & 0.0074 & 0.0077 & 0.0001 & 0.0256 \\
z8 & Gaussian Process & MSE & 0.0053 & 0.0055 & 0.001 & 0.0193 \\
z8 & Random Forest & MAE & 0.0507 & 0.0366 & 0.0056 & 0.1417 \\
z8 & Gaussian Process & MAE & 0.0454 & 0.0183 & 0.0244 & 0.087 \\
z1 & Random Forest & MSE & 0.0254 & 0.0407 & 0 & 0.1285 \\
z1 & Gaussian Process & MSE & 0.0973 & 0.178 & 0.0001 & 0.609 \\
z1 & Random Forest & MAE & 0.0715 & 0.093 & 0.0004 & 0.285 \\
z1 & Gaussian Process & MAE & 0.1405 & 0.1716 & 0.0071 & 0.5662 \\
\end{longtable}

\newpage{}

\section*{References}\label{references}
\addcontentsline{toc}{section}{References}

\protect\phantomsection\label{refs}

\end{document}